\documentclass[12pt]{article}
\usepackage{amssymb, amsmath, amsthm, amscd, graphicx, dsfont, mathabx}

\newtheorem{lem}{Lemma}[section]
\newtheorem{pro}{Proposition}[section]
\newtheorem{thm}{Theorem}[section]

\newtheorem{defn}{Definition}[section]
\font\toto=cmssdc10 scaled 1250

\def\B{\mathbb{B}}
\def\N {\mathbb{N}}
\def\R {\mathbb{R}}
\def\Z {\mathbb{Z}} 

\def\cA{{\mathcal A}}

\def\cE{{\mathcal E}} 
 
\def\cG {\mathcal{G}} 
\def\cL {\mathcal{L}}
\def\cP {\mathcal{P}}
\def\cR {\mathcal{R}}
\def\cV {\mathcal{V}} 
\def\cT{{\mathcal T}}
\def\cW{{\mathcal W}}

\def\tA{{\tt A}}
\def\tb{{\tt b}}

\def\colon{\,{:}\;}
\def\prf{\medbreak\noindent{\bf Proof.}\enspace}
\def\qed{\hspace*{\fill}\hbox{\vrule height 7pt \kern-.3pt
     \vbox{\hrule width 7pt
     \kern6.6pt\hrule width 7pt }\kern-.3pt\vrule height 7pt
     }\par}
\def\ra{\rightarrow}

\def\!{\mskip-\thinmuskip}

\begin{document}

\title{\bf Gibbs  measures on compact ultrametric spaces}

\author{C.-E. Pfister\footnote{E-mail: charles.pfister@epfl.ch}\\
Section of Mathematics,
Faculty of Basic Sciences, EPFL\\
CH-1015 Lausanne, Switzerland}


\maketitle
\noindent
{\bf Abstract:\,} 
One proves the equivalence of a Gibbs measure and a Gibbs conformal measure
for a dynamical system (G,X) when G is a countably infinite discrete group acting expansively on a
compact ultrametric space X.
As an application one proves for any beta-shift, that
the unique equilibrium measure for a function of  summable variation  is a Gibbs measure.

\vspace*{1.5cm}

\section{Introduction}\label{introduction}
\setcounter{equation}{0}

The notion of Gibbs measure is a basic notion both in Statistical Mechanics and Dynamical Systems. 
To put the results into perspective consider the
standard setting in Statistical Mechanics of lattice systems, which corresponds to the dynamical system $(G,X)$ where
$G=\Z^d$ and $X=\tA^{\Z^d}$ with $\tA$ a finite set. The action of the group $G$ on  $x\in X$, $k\mapsto x(k)$, is
by translations,
$$
(j,x)\mapsto T^jx\quad\text{where}\quad (T^jx)(k):=x(k+j)\,.
$$ 
Let $\Lambda_n=\{k\in G: \max\{|k_1|,\ldots,|k_d|\}\leq n\}$ and $\epsilon_n$, $n\geq 0$, be a strictly decreasing sequence, 
$\epsilon_0=1$ and $\lim_n\epsilon_n=0$. For the distance $d$ one chooses
$$
d(x,y)=\epsilon_m\quad\text{if}\;m=\min\{n\colon x(k)\not =y(k)\,\;\text{for some $k\in\Lambda_n$}\}\,.
$$
Then $(X,d)$ is ultrametric and the dynamical system $(G,X)$ is expansive, with expansive constant any $0<\delta<1$. The usual way to define Gibbs measures in this setting is to start with absolutely summable potentials $\Phi$, which are collections of $A$-local functions $\Phi_A$, $A\subset G$ finite, such that $\Phi_{T^k A}(x)=\Phi_A(T^k x)$ and $\sum_{A\ni 0}\|\Phi_A\|_\infty<\infty$ \cite{Ru}, \cite{Ge}. Associated to $\Phi$ is the continuous function
$$
f_\Phi=\sum_{A\ni 0}\frac{\Phi_A}{|A|}\,.
$$
In this paper one generalises this case in two directions.

Case (a). One assumes that $G$ is a countably infinite discrete set (no group structure)
and  $X$ is a closed subset of the product space $\prod_{g\in G}\tA_g$, where $\tA_g$ is a finite set for each $g\in G$. 
The potential $\Phi$ is replaced by a family of continuous functions 
$F=\{f_g\in C(X): g\in G\}$ such that for all finite subsets $\Lambda\subset G$
$$
\sup\big\{\sum_{g\not\in \Lambda}\big|f_g(x)-f_g(y)\big|\colon
x(g)=y(g)\,\;\forall g\not\in \Lambda\big\}<\infty\,.
$$
This corresponds to the setting for Gibbs measures in \cite{Ru}, but with the only assumption that $X$ is a closed subset of $\prod_{g\in G}\tA_g$.

Case (b). The metric space  $(X,d)$ is a compact ultrametric space and $G$ is a countably infinite discrete group. There is an action of $G$ on $X$ by homeomorphisms $x\mapsto g\cdot x$ and this action is expansive. 
The potential is replaced by a continuous function $f$ such that for all finite subsets $\Lambda\subset G$,
\begin{equation}\label{con}
\sup\big\{\sum_{g\not\in \Lambda}\big|f(g\cdot x)-f(g\cdot y)\big|\colon
d(g\cdot x,g\cdot y)\leq\varepsilon^*\,\;\forall g\not\in \Lambda\big\}<\infty\,,
\end{equation}
where $\varepsilon^*$ is an expansive constant.

There are two approaches to define Gibbs measures \cite{Ru}. Gibbs measures are defined by Dobrushin-Lanford-Ruelle equations \cite{Ge}. This is how they were first defined by Dobrushin \cite{Do} and by Lanford and Ruelle \cite{LR}. The second approach was initiated by Capocaccia \cite{Ca}.  Here the starting point is the notion of  
conjugate points. Conjugation defines an equivalence relation on $X$.
This leads to the notion of  Gibbs conformal measures.
The main result of this paper in section \ref{sectionmainthm} is that Gibbs conformal measures and Gibbs measures are equivalent if $X$ is a compact ultrametric space, theorem \ref{thmM5.2}.

In the second part of the paper, using the result of  section \ref{sectionmainthm} and  results of Climenhaga and Thompson \cite{CT}, 
one proves for the natural extension of any $\beta$-shift, when  $f$ is a continuous function verifying 
Bowen's condition and  condition \eqref{con} (with $G=\Z$ and the action by translation), that the unique equilibrium measure for $f$ is a Gibbs measure.

After putting this paper in {\it arXiv} I learned that Kimura proved already in 2015, for general subshifts over $\Z^d$, that a Gibbs conformal measure is a Gibbs measure (terminology of this paper), see theorem 5.30 \cite{Ki}. His proof is different from the proof given in subsection \ref{subsectionproof} part a). He also proved the converse, part b), but only for continuous cocycles, see theorem 5.32  \cite{Ki}. His results were extended to subshifts over a countable discrete group $G$ and theorem  \ref{thmM5.2} was proven in this context by L. Borsato and S. MacDonald. See theorems 5 and 6 in \cite{BM}. 

{\it Acknowledgement.\,} I thank R. Bissacot for  correspondence about Kimura's thesis.

\section{General setting}\label{setting}
\setcounter{equation}{0}

In the whole paper $G$ is a countably infinite discrete set. When $G$ is a group, the group operation is written $(g,h)\mapsto gh$ and $e$ is the neutral element of $G$. The set of all finite subsets of $G$ is denoted $\cP_f(G)$ and is directed by  $\supset$. For $\Lambda\subset G$, 
$\Lambda^c=G\setminus\Lambda$.
One considers two cases.
\begin{itemize}
\item[(a)]  Lattice system $X$.
For each $g\in G$ there is a finite set $\tA_g$ with the discrete distance $d$ on each $\tA_g$ and $X$ is a closed subset of the product space
$\prod_{g\in G}\tA_g$. An element $x\in X$ is $g\mapsto x(g)$, $x(g)\in\tA_g$.
\item[(b)] Dynamical system $(G,X)$. $G$ is a group and $(X,d)$ is a compact ultrametric space,
$$
d(x,z)\leq \max\{d(x,y),d(y,z)\}\quad\forall x,y,z\in X\,.
$$
The group $G$ acts on $X$ by homeomorphims and  the action of $G$ on $X$ is written
$(g,x)\mapsto g\cdot x$.
The action of $G$ is always expansive: there exists a positive constant $c$ such that
$$
x\not=y\implies \exists g\in G\,,\; d(g\cdot x,g\cdot y)>c\,.
$$
Without restricting the generality one assumes that $c<1$.
\end{itemize}
The set of real-valued continuous functions on $X$ with the sup-norm $\|\cdot\|_\infty$ is written $C(X)$. All measures on $X$ are Borel measures and
the set of probability measures is $M_1(X)$ with the weak topology. The group $G$ acts also on functions, $(g\cdot f)(x):=f(g\cdot x)$.

For $\Lambda\subset G$ set
\begin{eqnarray}\label{dLambda}
\text{case (a):} &\quad& d_\Lambda(x,y):=\sup\{d(x(g),y(g)):g\in \Lambda\}\\
\text{case (b):} &\quad& d_\Lambda(x,y):=\sup\{d(g\cdot x,g\cdot y):g\in \Lambda\}\,.
\end{eqnarray}
Let $0<\varepsilon^*\leq c$ be  a fixed  constant. 
For both cases
\begin{eqnarray}\label{WLambda}
[u]_\Lambda &:=& \{y\in X\colon d_\Lambda(u,y)\leq\varepsilon^*\}\nonumber\\
W_\Lambda^x &:=&\{y\in X\colon d_{\Lambda^c}(x,y)\leq\varepsilon^*\}\,.\nonumber
\end{eqnarray}
In case (a), since $c<1$, $[u]_\Lambda=\{y\colon y(g)=u(g)\,\;\forall g\in\Lambda\}$. In case (b)
the dependency of the notions and results on the choice of $\varepsilon^*$ is discussed in subsection
\ref {independence}. It is convenient to treat both cases together and  one omits $\varepsilon^*$ in the notations.
When $\Lambda\in\cP_f(G)$, $[u]_\Lambda$ is a \emph{cylinder set}.
The relation
$x\sim_\Lambda y \iff d_\Lambda(x,y)\leq\varepsilon^*$
is an equivalence relation. 
For finite $\Lambda$ 
$[u]_\Lambda$ is a clopen subset. Either $[u]_\Lambda$ and 
$[v]_\Lambda$ are disjoint or are equal. 
Given $x$ and $\Lambda$, 
\begin{equation}\label{equivWLambda}
y\in W_\Lambda^x\iff W_\Lambda^y=W_\Lambda^x\,.
\end{equation}
For case (a) this is obvious since $d(x(y),x(y))\leq\varepsilon^*$ means $x(g)=y(g)$. For case $(b)$ this is a consequence of the ultrametrity.

\begin{lem}\label{lemdistance}
Let $(G,X)$ be a dynamical system, case $(b)$. Let $\Lambda\in\cP_f(G)$ and
$$
\delta_\Lambda:=\sup\{d(x,y)\colon \text{$x,y\in X$, $d_\Lambda(x,y)\leq\varepsilon^*$}\}\,.
$$
i)  $\Lambda\supset\Lambda^\prime\implies \delta_{\Lambda}\leq \delta_{\Lambda^\prime}$.\\
ii) The net $\{\delta_\Lambda: \Lambda\in\cP_f(G)\}$ converges to $0$,
$\lim_{\Lambda\uparrow G}\delta_\Lambda=0$.\\
iii) Let $0<s<c$. Then 
$$
d(x,y)>s \implies \exists A_s\in\cP_f(G)\,,\, \exists g\in A_s\;\text{such that $d(g\cdot x,g\cdot y)>c$}\,.
$$
iv)  For $\Lambda\in\cP_f(G)$ there exists $\delta(\Lambda)>0$, such that 
$$
d(x,y)\leq\delta(\Lambda)\implies d_\Lambda(x,y)\leq\varepsilon^*\,.
$$
\end{lem}
\prf
Case i) is clear.\\
ii) Suppose that $\delta_\Lambda\geq \delta>0$ for all $\Lambda\in\cP_f(G)$. Let  $\{\Lambda_n\}$ be a sequence in $\cP_f(X)$
such that 
\begin{equation}\label{Lambdan}
e\in\Lambda_n\,,\quad\Lambda_n\subset\Lambda_{n+1}\,,\quad \bigcup_{n}\Lambda_n=G\,.
\end{equation}
Then there exist  
 $x_n,y_n$ such that
$$
d(x_n,y_n)\geq\delta/2\quad\text{and}\quad d(g\cdot x_n, g\cdot y_n)\leq\varepsilon^*\; \forall g\in\Lambda_n\,.
$$
Since $X$ is compact one may assume that $\lim_nx_n=x$ and $\lim_ny_n=y$. For any $g\in G$,
$$
d(g\cdot x, g\cdot y)=\lim_{n\ra \infty} d(g\cdot x_n, g\cdot y_n)\leq\varepsilon^*\quad\text{and}\quad d(x,y)=\lim_{n\ra\infty}d(x_n,y_n)\geq\delta/2\,,
$$
a contradiction with expansiveness. \\
iii) Suppose that $d(x,y)>s$ and  iii) is false. Let $\{\Lambda_n\}$ be as in \eqref{Lambdan}. Then for each $n$ there exist $x_n,y_n$ such that 
$$
d(x_n,y_n)\geq s\quad\text{and}\quad d(g\cdot x_n,g\cdot y_n)\leq c\quad \forall g\in \Lambda_n\,.
$$
One may assume that the sequences $x_n$ and $y_n$ converge  to $x$ and $y$ respectively. One has $d(x,y)\geq s$ and
$$
\forall n\;\text{and}\;\forall g\in \Lambda_n\quad
 d(g\cdot x_n,g\cdot y_n)\leq c\implies d(g\cdot x,g\cdot y)\leq c\quad\forall g\in G\,.
$$
This contradicts the fact that $c$ is an expansive constant.\\
iv) The homeomorphism $x\mapsto g\cdot x$  is uniformly continuous since $X$ is compact. Let
$$
\delta_g:=\sup\{\delta: d(x,y)\leq\delta\implies d(g\cdot x,g\cdot y)\leq\varepsilon^*\}\,.
$$
Continuity of the map $x\mapsto g\cdot x$ implies $\delta_g>0$. Since $\Lambda$ is finite,
$$
\delta(\Lambda):=\min\{\delta_g\colon g\in\Lambda\}>0\,.
$$
\qed

Two points $x,y\in X$ are \emph{conjugate} if for any $\varepsilon>0$ there exists a finite subset
$\Lambda_\varepsilon$ such that $d_{\Lambda^c}(x,y)\leq\varepsilon$ for all $\Lambda\supset\Lambda_\varepsilon$.
The set of conjugate points to $x$ is written
$$
W(x):=\{z\in X\colon \lim_{\Lambda\uparrow G}d_{\Lambda^c}(x,z)=0\}\,.
$$
Conjugacy  defines an  equivalence relation  $\sim$, called \emph{tail equivalence relation} or \emph{homoclinic relation} \cite{Me},
$$
\cR(X)=\big\{(x,y)\in X\times X:
\lim_{\Lambda\uparrow G}d_{\Lambda^c}(x,y)=0\big\}\,.
$$
The equivalence class of $x$ is $[x]_\sim=W(x)$. 

\begin{lem}\label{lemconjugacy}
$$
[x]_\sim=\bigcup_{\Lambda \in\cP_f(G)}W_\Lambda^x\,.
$$
The sets $W_\Lambda^x$ are finite and each equivalence class  $[x]_\sim$ is countable. 
\end{lem}

\prf
The case (a)  is clear since $d_{\Lambda^c}(x,y)\leq\varepsilon^*$ implies 
$x(g)=y(g)$ for all $g\not\in \Lambda$. \\
Case (b).
For $\Lambda$ finite there is a finite partition of $X$ into subsets $[u]_\Lambda$, say
$[u_i]$, $i=1,\ldots N$. Let $z\in W_\Lambda^x$; then
$z\in [u_i]_\Lambda$ for some $i$. The set $[u_i]_\Lambda\cap W_\Lambda^x$ has cardinality at most $1$, 
since $y,z\in [u_i]_\Lambda\cap W_\Lambda^x$ implies $d(g\cdot y,g\cdot z)\leq \varepsilon^*$ for all $g\in G$.
Hence $|W_\Lambda^x|\leq N$.

Let $\{\Lambda_n\}\subset\cP_f(G)$ be as in \eqref{Lambdan} and
$y\in \bigcup_{\Lambda \in\cP_f(G)}W_\Lambda^x$.  Suppose that $y\in W_{\Lambda^\prime}^x$. Hence
$y\in W_{\Lambda_n}^x$ for all $\Lambda_n\supset\Lambda^\prime$.
 Let $\varepsilon>0$ be given and $\Lambda_m$ such that 
$\delta_{\Lambda_m}\leq\varepsilon$ (lemma \ref{lemdistance} ii). Let
$$
A_\varepsilon:=\{g\in G\colon  \Lambda_mg\cap\Lambda^\prime=\emptyset\}\,.
$$
The complement of $A_\varepsilon$ is finite. Indeed, if $\Lambda_mg\cap\Lambda^\prime\not=\emptyset$, then there exist
$k\in\Lambda^\prime$ and  $\ell\in\Lambda_m$ such that $k=\ell g$, or $g=\ell^{-1}k$. Therefore
$$
A_\varepsilon^c\subset \bigcup_{k\in\Lambda^\prime}\{g=\ell^{-1}k: \ell\in\Lambda_m\}\,.
$$
For $n$ large enough  $A^c_\varepsilon\subset\Lambda_n$.
Hence for all $g\in \Lambda^c_n$, 
$$
\sup_{h\in\Lambda_m}d(h\cdot(g\cdot x),h\cdot(g\cdot y))\leq\varepsilon^*\implies d(g\cdot x,g\cdot y)\leq\delta_{\Lambda_m}\leq\varepsilon\,.
$$
Thus $d_{\Lambda^c_n}(x,y)\leq\varepsilon$, proving that $x\sim y$.

Conversely, if $x\sim z$, then there exists $\Lambda^\prime\in\cP_f(G)$ such that $d(g\cdot x,g\cdot z)\leq\varepsilon^*$ for all $g\not\in\Lambda^\prime$, that is $z\in W_{\Lambda^\prime}^x$. 
\qed

By definition, a set $A\subset X$ is \emph{$\Lambda$-saturated} if
$$
x\in A\implies W_\Lambda^x\subset A\quad\text{i.e.}\quad A=\bigcup_{x\in A}W_\Lambda^x\,.
$$
Complements and unions of $\Lambda$-saturated sets are $\Lambda$-saturated.

\begin{defn}\label{defnBspace}
Let $f\in C(X)$ and set
\begin{equation}\label{eq7} 
\Delta_{\Lambda}(f):=\sup\big\{\sum_{g\not\in \Lambda}\big|f(g\cdot x)-f(g\cdot y)\big|\colon x,y\in X\,,\;
\text{$d_{\Lambda^c}(x, y)\leq\varepsilon^*$}\big\}\,.
\end{equation}
Then
$$
\B(X):=\{f\in C(X)\colon \Delta_{\Lambda}(f)<\infty\,\;\text{for all finite $\Lambda$}\}\,.
$$
\end{defn}
For case (a), instead of $f\in\B(X)$ one has a set $F:=\{f_g\in C(X):g\in G\}$ such that for all $\Lambda\in\cP_f(G)$
\begin{equation}\label{eq7b} 
\Delta_{\Lambda}(F):=\sup\big\{\sum_{g\not\in \Lambda}\big|f_g(x)-f_g(y)\big|\colon
\text{$x(g)=y(g)$ $\forall g\in \Lambda^c$}\big\}<\infty\,.
\end{equation}

\section{Gibbs measures}\label{sectionGibbs}
\setcounter{equation}{0}

Let $f\in\B(X)$ and   $\Lambda\in\cP_f(G)$. Define
\begin{eqnarray}\label{eqprobabilitykernelnew}
\pi_\Lambda^f(y|x)
&:=&
\begin{cases}
\displaystyle \frac{\exp\sum_{g\in G}\big(f(g\cdot y)-f(g\cdot x)\big)}
{\sum_{y^\prime\in W^x_\Lambda}\exp\sum_{g\in G}\big(f(g\cdot y^\prime)-f(g\cdot x)\big)} &\text{if $y\in W_\Lambda^x$}\\
0 &\text{if $y\not\in W_\Lambda^x$.}
\end{cases}
\end{eqnarray}
This expression is well-defined because  $\sum_{g\in G}\big(f(g\cdot y)-f(g\cdot x)\big)$ is absolutely summable; it  is a probability measure on $W_\Lambda^x$. Moreover
\begin{equation}\label{eqequality}
\pi_\Lambda^f(y|x)=\pi_\Lambda^f(y|x^\prime)\quad\text{if $x^\prime\in W_\Lambda^x$}\,.
\end{equation}
Indeed, $x^\prime\in W_\Lambda^x$ if and only if $W_\Lambda^x=W_\Lambda^{x^\prime}$ and
$$
\varphi(g\cdot y)-\varphi(g\cdot x)=\big(\varphi(g\cdot y)-\varphi(g\cdot x^\prime)\big)+
\big(\varphi(g\cdot x^\prime)-\varphi(g\cdot x)\big)\,.
$$
A convenient expression for $\pi_\Lambda^f(y|x)$, when $y\in W_\Lambda^x$, is
\begin{equation}\label{eqconvenient}
\pi_\Lambda^f(y|x)=
\displaystyle\frac{1}{\sum_{y^\prime\in W^x_\Lambda}\exp\sum_{g\in G}\big(f(g\cdot y^\prime)-f(g\cdot y)\big)}\,.
\end{equation}
When $X$ is a lattice system, instead of  $\pi_\Lambda^f(x|y)$ one has $\pi_\Lambda^F(x|y)$, $F\subset C(X)$ verifying \eqref{eq7b}, and 
$\pi_\Lambda^F(x|y)$ is obtained by replacing in \eqref{eqprobabilitykernelnew} $f(g\cdot x)$ by $f_g(x)$.  For the exposition one considers case (b).

Let $\varphi$ be a bounded  function on $X$ and set
\begin{equation}\label{MLambda}
M_\Lambda(\varphi)(x):=
\sum_{y\in W_\Lambda^x}\varphi(y)\pi_\Lambda^f(y|x)\,.
\end{equation}
The function $M_\Lambda(\varphi)$  is constant on $W_{\Lambda}^x$.

\begin{pro}\label{proprojection}
Let $\Lambda\subset\Lambda^\prime$. Then
$$
M_{\Lambda^\prime}(\varphi)=M_{\Lambda^\prime}(M_{\Lambda}(\varphi))\,.
$$
\end{pro}

\prf
Let $\Lambda\subset\Lambda^\prime$.  Set
$$
Z_{\Lambda^\prime}(x):=\sum_{z\in W_{\Lambda^\prime}^x}\exp\sum_{g\in G}(f(g\cdot z)-f(g\cdot x))\,.
$$
 By ultrametricty for each $y\in W^x_{\Lambda^\prime}$, $W_\Lambda^y\subset W^x_{\Lambda^\prime}$.
The set $W^x_{\Lambda^\prime}$ is partitioned according to the equivalence relation: $y,z\in W^x_{\Lambda^\prime}$ are equivalent iff
$W_\Lambda^y=W_\Lambda^z$. For each equivalence class $C_i$ one chooses one element $\bar y_i$. Relation \eqref{equivWLambda} implies that
$C_i=W^{\bar y_i}_{\Lambda}$.
\begin{eqnarray*}
M_\Lambda(\varphi)(\bar y_i)
&=&
\displaystyle\frac{\sum_{y\in W_\Lambda^{\bar y_i}}\varphi(y)\exp\sum_{g\in G}(f(g\cdot y)-f(g\cdot x))}
{\sum_{z\in W_\Lambda^{\bar y_i}}\exp\sum_{g\in G}(f(g\cdot z)-f(g\cdot x))} \\
& =&\displaystyle\frac{\sum_{y\in C_i}\varphi(y)\exp\sum_{g\in G}(f(g\cdot y)-f(g\cdot \bar y_i))}
{\sum_{z\in C_i}\exp\sum_{g\in G}(f(g\cdot z)-f(g\cdot \bar y_i))}\,.
\end{eqnarray*}
For all $z\in C_i$, $M_\Lambda(\varphi)(\bar y_i)=M_\Lambda(\varphi)(z)$.
Therefore, setting
$$
Z_\Lambda(\bar y_i,x)=\sum_{z\in C_i}\exp\sum_{g\in G}(f(g\cdot z)-f(g\cdot x))\,,
$$
\begin{eqnarray*}
Z_{\Lambda^\prime}(x) M_{\Lambda^\prime}(\varphi)(x)
&=&
Z_{\Lambda^\prime}(x)\sum_i\sum_{y\in C_i}\varphi(y)\pi_{\Lambda^\prime}^f(y|x)\\
&=&
\sum_i Z_\Lambda(\bar y_i,x)
\displaystyle\frac{\sum_{y\in C_i}\varphi(y)\exp\sum_{g\in G}(f(g\cdot y)-f(g\cdot x))}{Z_\Lambda(\bar y_i,x)}\\
&=&
\sum_i Z_\Lambda(\bar y_i,x) M_\Lambda(\varphi)(\bar y_i)\\
&=&
\sum_i\Big(\sum_{z\in C_i}\exp\sum_{g\in G}(f(g\cdot z)-f(g\cdot x))\Big)M_\Lambda(\varphi)(z)\\
&=&
Z_{\Lambda^\prime}(x) M_{\Lambda^\prime}(M_\Lambda(\varphi))(x)\,.
\end{eqnarray*}
\qed

\begin{defn}\label{defnweakdependence}
The dynamical system $(G,X)$, or the lattice system $X$, has the \emph{weak dependence property} if for all 
$\Lambda\in\cP_f(G)$ there exists  $\bar \Lambda\in\cP_f(G)$ such that
$\bar \Lambda\supset\Lambda$ and for all $v\in X$
\begin{eqnarray}\label{eqdependence}
d_{\bar \Lambda\setminus\Lambda}(x,x^\prime)\leq\varepsilon^*
&\implies&
\big([v]_{\Lambda}\cap W_{\Lambda}^x\not=\emptyset\iff [v]_{\Lambda}\cap W_{\Lambda}^{x^\prime}\not=\emptyset\big)\,.
\end{eqnarray}
\end{defn}

If $X$ has the weak dependence property and $d_{\bar \Lambda\setminus\Lambda}(x,x^\prime)\leq\varepsilon^*$, then there is a natural bijection  $W_\Lambda^x\ra W_{\Lambda}^{x^\prime}$:
\begin{equation}\label{bijection}
z\in [w]_\Lambda\cap W_\Lambda^x\mapsto  z^\prime\in [w]_\Lambda\cap W_\Lambda^{x^\prime}\,.
\end{equation}
The  \emph{topological Markov property} of \cite{BGMT} for  the dynamical system $(G,X)$ or the lattice system $X$
reads as follows:
for any  $\Lambda\in \cP_f(G)$ there exists a  subset $\bar\Lambda\in\cP_f(G)$ such that 
$\bar\Lambda\supset\Lambda$ and for $u,v,y\in X$,
if $d_{\bar\Lambda\setminus\Lambda}(u,v)\leq\varepsilon^*$, then
\begin{equation*}
\exists z: d_{\bar\Lambda}(u,z)\leq\varepsilon^*\;\text{and}\; d_{\bar\Lambda^c}(y,z)\leq\varepsilon^*
\iff
\exists  z^\prime: d_{\bar\Lambda}(v,z^\prime)\leq\varepsilon^*\;\text{and}\; d_{\bar\Lambda^c}(y,z^\prime)\leq\varepsilon^*\,.
\end{equation*}
 
\begin{lem}\label{equivalence}
The weak dependence property and the topological Markov property are equivalent.
\end{lem}

\prf
Suppose that $(G,X)$ has the weak dependence property and that $d_{\bar\Lambda\setminus\Lambda}(u,v)\leq\varepsilon^*$.
Let $y\in X$ and suppose that there exists  $z$, $d_{\bar\Lambda}(u,z)\leq\varepsilon^*$ and $d_{\bar\Lambda^c}(y,z)\leq\varepsilon^*$.
Then $d_{\bar\Lambda\setminus\Lambda}(v,z)\leq\varepsilon^*$ and
$$
v\in [v]_\Lambda\cap W_\Lambda^v \implies [v]_\Lambda\cap W_\Lambda^z\not=\emptyset\,.
$$
Hence, there exists $z^\prime$ such that $d_{\Lambda}(v,z^\prime)\leq\varepsilon^*$ and $d_{\Lambda^c}(z,z^\prime)\leq\varepsilon^*$.
It follows that    $d_{\bar\Lambda}(v,z^\prime)\leq\varepsilon^*$ and $d_{\bar\Lambda^c}(y,z^\prime)\leq\varepsilon^*$. 

Suppose that $(G,X)$ has the topological Markov property. Let $x,x^\prime$, $d_{\bar\Lambda\setminus \Lambda}(x,x^\prime)\leq\varepsilon^*$ and
$[v]_\Lambda\cap W_\Lambda^x\ni z$. 
It follows that
 $d_{\bar\Lambda\setminus\Lambda}(x^\prime,z)\leq\varepsilon^*$. Applying the Markov property with $u=x^\prime$ and 
 $v=z$, and taking $y=x^\prime$, one concludes the existence of $z^\prime$,
 $$
d_{\bar\Lambda}(v,z^\prime)=d_{\bar\Lambda}(z,z^\prime)\leq\varepsilon^*\;\text{and}\; 
d_{\bar\Lambda^c}(y,z^\prime)=d_{\bar\Lambda^c}(x^\prime,z^\prime)\leq\varepsilon^*\,.
 $$
Hence $d_{\Lambda}(v,z^\prime)\leq\varepsilon^*$ and $d_{\Lambda^c}(x^\prime,z^\prime)\leq\varepsilon^*$, that is
$[v]_\Lambda\cap W_\Lambda^{x^\prime}\not=\emptyset$.
\qed

\begin{pro}\label{procontinuityMLambda}
If $(X,G)$ or $X$ has the weak dependence property, then $M_\Lambda$ maps $C(X)$ into $C(X)$. 
\end{pro}

\prf
Case (b). 
Let $\varphi\in C(X)$, $\Lambda\in\cP_f(G)$ and $\varepsilon>0$. One may assume that $\varphi\geq 0$.
 Choose $\Lambda^\prime\in\cP_f(G)$  such that 
$\Lambda^\prime\supset\bar\Lambda$ and
$$
\sup\big\{\sum_{g\not\in \Lambda^\prime}\big|f(g\cdot x)-f(g\cdot y)\big|\colon
d_{\Lambda^c}(x,y)\leq \varepsilon^*\big\}\leq\varepsilon\,.
$$
Each continuous function on $X$ is uniformly continuous. Let $\varepsilon^\prime=|\Lambda^\prime|^{-1}\varepsilon$ and choose $\delta>0$ so that 
$$
d(x,y)\leq\delta\implies |\varphi(x)-\varphi(y)|\leq\varepsilon \quad \text{and}\quad 
|f(g\cdot x)-f(g\cdot y)|\leq\varepsilon^\prime\,\;\forall g\in\Lambda^\prime\,.
$$
Choose  in $\cP_f(G)$ $\Lambda^{\prime\prime}\supset\bar\Lambda$ so that 
$d_{\Lambda^{\prime\prime}}(x,y)\leq \varepsilon^*$ implies $d(x,y)\leq\delta$ (lemma \ref{lemdistance} ii).
The weak dependence property implies that there is a bijection $\theta: W_\Lambda^x\ra W_\Lambda^y$ such that
$$
z\in W_\Lambda^x\cap[u]_\Lambda\implies \theta(z)\in W_\Lambda^y\cap[u]_\Lambda\,.
$$
Therefore $d_{\Lambda}(z,\theta(z))\leq\varepsilon^*$ and 
$$
d_{\Lambda^{\prime\prime}\setminus\Lambda}(z,\theta(z))\leq
\max\{d_{\Lambda^{\prime\prime}\setminus\Lambda}(z,x),d_{\Lambda^{\prime\prime}}(x,y),
d_{\Lambda^{\prime\prime}\setminus\Lambda}(y,,\theta(z)\}\leq\varepsilon^*\,,
$$
which implies that $d_{\Lambda^{\prime\prime}}(z,\theta(z))\leq\varepsilon^*$ and $d(z,\theta(z))\leq\delta$. For $d_{\Lambda^{\prime\prime}}(x,y)\leq \varepsilon^*$
\begin{eqnarray*}
& &\big|\sum_{g\in G}\big(f(g\cdot z)-f(g\cdot x)\big)-
\sum_{g\in G}\big(f(g\cdot \theta(z))-f(g\cdot y)\big)\big|\\
& \leq&
\big|\sum_{g\in \Lambda^\prime}\big(f(g\cdot z)-f(g\cdot\theta(z))\big)\big|
+
\big|\sum_{g\in \Lambda^\prime}\big(f(g\cdot x)-f(g\cdot y)\big)\big|+2\varepsilon\\
&\leq&
2|\Lambda^\prime|\varepsilon^\prime+2\varepsilon=4\varepsilon\,.
\end{eqnarray*}
From these estimates and  $d(x,y)\leq \delta$  one obtains 
\begin{eqnarray*}
M_\Lambda(\varphi)(x)&=&\sum_{z\in W_\Lambda^x}\pi^f_\Lambda(z|x)\big(\varphi(z)-\varphi(\theta(z))\big)+\sum_{z\in W_\Lambda^x}\pi^f_\Lambda(z|x)\varphi(\theta(z))\\
&\leq&
\varepsilon+ \sum_{z\in W_\Lambda^x}\pi^f_\Lambda(z|x)\varphi(\theta(z))\\
&\leq&
\varepsilon+{\rm e}^{8\varepsilon}\sum_{z\in W_\Lambda^x}\pi^f_\Lambda(\theta(z)|y)\varphi(\theta(z))\\
&=&
\varepsilon+{\rm e}^{8\varepsilon}\sum_{z^\prime\in W_\Lambda^y}\pi^f_\Lambda(z^\prime|y)\varphi(z^\prime)
=\varepsilon+{\rm e}^{8\varepsilon}M_\Lambda(\varphi)(y)\,.
\end{eqnarray*}
Under the same conditions, 
$ M_\Lambda(\varphi)(x)\geq -\varepsilon+{\rm e}^{-8\varepsilon}M_\Lambda(\varphi)(y)$.
\qed

\begin{defn}\label{defnGibbs}
Let $f\in\B(X)$ or $F\subset C(X)$ verifying \eqref{eq7b}.
A \emph{Gibbs measure}  is a probability measure $\mu$ such that
$$
\int \varphi\,d\mu=\int M_\Lambda(\varphi)\,d\mu\quad\forall \text{$\Lambda\in\cP_f(G)$  and all $\varphi\in C(X)$.}
$$
\end{defn}

\begin{thm}\label{thmexistence}
If $(G,X)$ or $X$ has the weak dependence property, then there exists a Gibbs measure.
\end{thm}

\prf
This follows from propositions \ref{proprojection} and \ref{procontinuityMLambda} by a standard compacity argument.
\qed

\section{Gibbs conformal measures}\label{sectionconformal}
\setcounter{equation}{0}

The reference for section \ref{sectionconformal} is \cite{Me}.
Let $f\in\B(X)$
 and define
\begin{equation}\label{eqdefcocycle}
\psi_f:\cR(X)\ra\R\,,\;(x,y)\mapsto \psi_f(x,y):=\sum_{g\in G}\big(f(g\cdot y)-f(g\cdot x)\big)\,.
\end{equation}
The map $\psi_f$ is well-defined on $\cR(X)$:  by definition of $\sim$, if $n$ is large enough,
$$
x\sim y \implies d(g\cdot x,g\cdot y)\leq\varepsilon^*\,\;\forall g\in\Lambda^c_n\,,
$$
so that
 the sum in \eqref{eqdefcocycle} is absolutely summable.
The map $\psi_f$ is a cocycle,
$$
\psi_f(x,y)+\psi_f(y,z)=\psi_f(x,z)\,.
$$
In the case of the lattice let $F\subset C(X)$ verifying \eqref{eq7b}, and instead of $\psi_f$ define
\begin{equation}\label{eqdefcocycleb}
\psi_F:\cR(X)\ra\R\,,\;(x,y)\mapsto \psi_f(x,y):=\sum_{g\in G}\big(f_g(y)-f_g(x)\big)\,.
\end{equation}
Let $\Lambda\in\cP_f(G)$ and $u,v\in X$. 
Set
\begin{equation}\label{eqinvolutionnew}
\phi_{\Lambda}^{u,v}(x):=\begin{cases}
y &\text{if $x\in [u]_{\Lambda}$ and  $y\in [v]_{\Lambda}\cap W_{\Lambda}^x$}\\
z &\text{if $x\in [v]_{\Lambda}$ and  $z\in [u]_{\Lambda}\cap W_{\Lambda}^x$}\\
x&\text{otherwise.}
\end{cases}
\end{equation}
The map $\phi_{\Lambda}^{u,v}$ is well-defined since $|[v]_{\Lambda}\cap W_{\Lambda}^x|\leq 1$ and 
$|[u]_{\Lambda}\cap W_{\Lambda}^x|\leq 1$. 
This class of transformations generate the tail-equivalence relation $\sim$,
$$
x\sim y\iff \exists\,\phi_{\Lambda}^{u,v}\,,\; y=\phi_{\Lambda}^{u,v}(x)\,.
$$
Indeed, if $y\sim x$, then $y\in W_\Lambda^x$ for some $\Lambda$. Hence, there exist $u$ and $v$ such that 
$x\in[u]_\Lambda$ and $y\in[v]_\Lambda\cap W_\Lambda^x$, i.e.  $y=\phi_{\Lambda}^{u,v}(x)$. If $y=\phi_{\Lambda}^{u,v}(x)$, then
$x\sim y$ (lemma \ref{lemconjugacy}).

\begin{lem}\label{leminvolutionnew}
The map $\phi_{\Lambda}^{u,v}$ is an involution.
\end{lem}

\prf
Either $x$ is a fixed point of  $\phi_{\Lambda}^{u,v}$ or there exist $u\not=v$ such that
$x\in  [u]_{\Lambda}$ and  $y=\phi_{\Lambda}^{u,v}(x)\not=x$. Claim $\phi_{\Lambda}^{u,v}(y)=x$. 
By definition of $\phi_{\Lambda}^{u,v}$
$\{y\}=[v]_{\Lambda}\cap W_{\Lambda}^x$. Thus $y\in[v]_\Lambda$ and 
$$
x\in[u]_\Lambda\quad\text{and}\quad d(g\cdot x, g\cdot y)\leq\varepsilon^*\,\;\forall g\in \Lambda^c
\implies x\in [u]_\Lambda\cap W_\Lambda^y\,.
$$
\qed

\begin{defn}\label{defnGibbs-conformal}
Let $f\in\B(X)$.  
A probability measure $\mu$ is  $(\psi_f,\cR)$-conformal if for any $\phi=\phi_\Lambda^{u,v}$
$$
\frac{d\mu\circ\phi}{d\mu}(x)=\exp \psi_f(x,\phi(x))\quad\text{$\mu$-a.s.}\,.
$$
The cocycle $\psi_f$ is a \emph{Gibbs cocycle} and $\mu$ is a \emph{Gibbs conformal  measure}. 
In the case of the lattice $X$, $F\subset C(X)$ verifies \eqref{eq7b} and  the cocycle is $\psi_F$.
\end{defn}

In general the involutions $\phi_\Lambda^{u,v}$ are not continuous maps. The next proposition gives a sufficient condition for
the continuity of the maps  $\phi_\Lambda^{u,v}$.

\begin{pro}\label{procontinuity}
If  the weak dependence property holds, then the  involutions $\phi_{\Lambda}^{u,v}$ are homeomorphisms.
\end{pro}

\prf
Given $\varepsilon>0$, let $\Lambda_\varepsilon\in\cP_f(G)$ be such that 
$\delta_{\Lambda_\varepsilon}\leq\varepsilon$ (lemma \ref{lemdistance} ii).
Let $x_n\ra x$. Choose $\Lambda^{\prime\prime}\in\cP_f(G)$ so that
$\Lambda^{\prime\prime}\supset\bar\Lambda\supset\Lambda$ and $\Lambda^{\prime\prime}\supset \Lambda_\varepsilon$. 
For $n$ large enough
$d(x_n,x)\leq\delta(\Lambda^{\prime\prime})$, so that  (lemma \ref{lemdistance} iv)
$$
d(g\cdot x_n,g\cdot x)\leq\varepsilon^*\,\;\forall g\in\Lambda^{\prime\prime}\,.
$$
Set $y=\phi_{\Lambda}^{u,v}(x)$ and $y_n=\phi_{\Lambda}^{u,v}(x_n)$. If $x,x_n\not\in [u]_\Lambda$ or $x,x_n\not\in [v]_\Lambda$, then
$y_n=x_n$ and $y=x$, so that $\lim_n\phi_{\Lambda}^{u,v}(x_n)=\phi_{\Lambda}^{u,v}(x)$.
Suppose that $x,x_n\in[u]_\Lambda$. By the weak dependence property, either $[v]_\Lambda\cap W_\Lambda^x=\emptyset$ and $[v]_\Lambda\cap W_\Lambda^{x_n}=\emptyset$, in which case  $y_n=x_n$ and $y=x$, or 
$y\in [v]_\Lambda\cap W_\Lambda^x$ and $y_n\in [v]_\Lambda\cap W_\Lambda^{x_n}$, so that 
$$
d(g\cdot y_n,g\cdot y)\leq \varepsilon^*\,\;\forall g\in \Lambda^{\prime\prime}\implies d(y_n,y)\leq\varepsilon\,.
$$
Since  $\phi_{\Lambda}^{u,v}$ is an involution, it is also an homeomorphism.
\qed

\section{Equivalence of  Gibbs conformal and Gibbs measures}\label{sectionmainthm}
\setcounter{equation}{0}

The main theorem is
 
\begin{thm}\label{thmM5.2}
Let $f\in\B(X)$ or $F\subset C(X)$ verify \eqref{eq7b}. Then
a probability measure is a Gibbs measure if and only if it is 
a Gibbs conformal measure.
\end{thm}

Let $\Lambda\in\cP_f(G)$.
There exists a unique finite partition $\cA_\Lambda$ of $X$ into subsets $[u]_\Lambda$.
For  each $[u]_\Lambda\in\cA_{\Lambda}$
one chooses an element, so that
$\cA_\Lambda=\{[u_1]_\Lambda,[u_2]_\Lambda,\ldots\}$.
Set
$$
\cW_\Lambda(x):=\{ [u_i]_\Lambda\in\cA_\Lambda\colon [u_i]_\Lambda\cap W_\Lambda^x\not=\emptyset\}\,.
$$
For a given $\cW\subset\cA_\Lambda$ set
$$
A_{\cW}:=\{x\in X\colon \cW_\Lambda(x)=\cW\}\,.
$$
The set $A_{\cW}$ is $\Lambda$-saturated.

\begin{lem}\label{lemmeasurable}
The sets $A_{\cW}$ define a measurable partition of $X$ and
the maps $\phi_\Lambda^{u,v}$ are measurable. If $\{[u]_\Lambda, [v]_\Lambda\}\subset\cW$, then for all $\Lambda$-saturated sets $B\subset A_{\cW}$, 
$\phi_\Lambda^{u,v}(B)=B$.
\end{lem}

\prf
i) $A_{\cW}$ is measurable.
Let $\{B_\ell\}$ be an increasing sequence of finite subsets of $\Lambda^c$ such that 
$\bigcup_\ell B_\ell=\Lambda^c$. For $[u_i]_\Lambda\in\cA_\Lambda$ set
$$
A_\ell(u_i):=\bigcup_{\substack{v:\\ [v]_{B_\ell}\cap[u_i]_\Lambda\not=\emptyset}}[v]_{B_\ell}\quad\text{and} \quad 
A(u_i):=\bigcap_\ell A_\ell(u_i)\,.
$$
$A(u_i)$ is closed ($A_\ell(u_i)$ is a finite union of clopen subsets); if $x\in A(u_i)$, then for all $\ell$, $x\in [v]_{B_\ell}$ for some $v$,
and for that $v$ there exists $x_\ell\in [v]_{B_\ell}\cap[u_i]_\Lambda$.
Ultrametricity implies  that $[x]_{B_\ell}=[v]_{B_\ell}$, so that $x_\ell\in [x]_{B_\ell}\cap[u_i]_\Lambda$. Since the sequence 
$\{B_\ell\}$ is increasing, by compactness there exists a subsequence
$\{x_{\ell_i}\}$ converging to $y$ and $y\in [x]_{B_\ell}\cap[u_i]_\Lambda$ for all $\ell$. Hence
$d_\Lambda(y, u_i)\leq\varepsilon^*$ and  $d_{\Lambda^c}( y,x)\leq\varepsilon^*$, i.e. $W_\Lambda^x\cap[u_i]_\Lambda=\{y\}$, and
$$
A(u_i)\subset \{x\in X: W_\Lambda^x\cap[u_i]_\Lambda\not=\emptyset\}\,.
$$ 
 Conversely, if $y\in W_\Lambda^x\cap[u_i]_\Lambda$, then
$y\in A(u_i)$, so that $A(u_i)=\{x: W_\Lambda^x\cap[u_i]_\Lambda\not=\emptyset\}$.
For $\cW\subset\cA_\Lambda$,
$$
A_\cW=\Big(\bigcap_{[u_i]_\Lambda\in\cW} A(u_i)\Big)\setminus\Big(\bigcup_{[u_j]_\Lambda\not\in\cW}A(u_j)\Big)\,,
$$
proving the measurability of $A_\cW$.

\noindent
ii) $\phi_\Lambda^{u,v}$ is measurable.
Let $E$ be a closed subset and $E^s=\bigcup_{x\in E}W_\Lambda^x$ the $\Lambda$-saturated set generated by $E$.
Let
$$
\bar E_\ell:=\bigcup_{x\in E}[x]_{B_\ell}\quad\text{and}\quad \bar E:=\bigcap_\ell \bar E_\ell\,.
$$
$\bar E$ is measurable ($\bar E_\ell$ is a finite union of cylinder sets) and $\bar E=E^s$. Indeed, suppose that $y\in E^s$ i.e. $y\in W_\Lambda^x$ for some $x\in E$.
Then $y\in [x]_{B_\ell}$ for all $\ell$, so that
$y\in\bar E$.  
Conversely, if $z\in \bar E$, then for any $\ell$ there exists $x_\ell\in E$ such that $z\in[x_\ell]_{B_\ell}$, or equivalently by ultrametricity, $x_\ell\in[z]_{B_\ell}$. By compactness there exists a subsequence $\{x_{\ell_i}\}$ converging to some $x\in E$ since $E$ is closed.
Therefore $z\in W_\Lambda^x$ i.e. $z\in E^s$.

Recall that $\phi_\Lambda^{u,v}(x)=x$ if $x\in[u]_{\Lambda}$ and  $[v]_\Lambda\cap W_\Lambda^x=\emptyset$
or  $x\in[v]_{\Lambda}$ and  $[u]_\Lambda\cap W_\Lambda^x=\emptyset$. The set $E$ is decomposed into
$$
E_1:=E\cap\bigcup\big\{A_\cW\colon \cW\not\supset\{[u]_\Lambda,[v]_\Lambda\}\big\}
$$
and
$$
E_2:=E\cap \bigcup\big\{A_\cW\colon \cW\supset\{[u]_\Lambda,[v]_\Lambda\}\big\}\,.
$$
The map $\phi_\Lambda^{u,v}$ is the identity map on $E_1$ and on $E_2\cap([u]_\Lambda\cup[v]_\Lambda)^c$.
On the other hand
$$
E_2^s=E^s\cap\bigcup\big\{A_\cW\colon \cW\supset\{[u]_\Lambda,[v]_\Lambda\}\big\}\,,
$$
and
\begin{eqnarray*}
(\phi_\Lambda^{u,v})^{-1}(E_2\cap[u]_\Lambda )
&=&
\begin{cases}
E_2^s\cap[v]_\Lambda &\text{if $E_2\cap[u]_\Lambda \not=\emptyset$}\\
\emptyset &\text{if $E_2\cap[u]_\Lambda =\emptyset$}
\end{cases}\\
(\phi_\Lambda^{u,v})^{-1}(E_2\cap[v]_\Lambda )
&=&
\begin{cases}
E_2^s\cap[u]_\Lambda &\text{if $E_2\cap[v]_\Lambda \not=\emptyset$}\\
\emptyset &\text{if $E_2\cap[v]_\Lambda =\emptyset$.}
\end{cases}
\end{eqnarray*}
Therefore, for all closed subsets $E$, $(\phi_\Lambda^{u,v})^{-1}(E)$ is measurable.\\
iii) Finally,  either $[u]_\Lambda\cap W_\Lambda^x=\emptyset$, or
there exists a unique $x_u\in [u]_\Lambda\cap W_\Lambda^x$.
Since $\phi_\Lambda^{u,v}$ permutes $x_u\in B$ and $x_v\in B$, $\phi_\Lambda^{u,v}(B)=B$.
\qed

\subsection{Proof of theorem \ref{thmM5.2}}\label{subsectionproof}
a) Let  $\mu$ be a Gibbs conformal measure with cocycle $\psi_f$. 
Fix $\cW\ni[u]_\Lambda$ and let $B$ be a measurable $\Lambda$-saturated subset of $A_\cW$. 
If  $[v]_\Lambda\in\cW$, then
$\phi_\Lambda^{u,v}([u]_\Lambda\cap B)=[v]_\Lambda\cap B$ and 
\begin{eqnarray}
\mu\circ\phi^{u,v}([u]_\Lambda\cap B)&=&\int_{\phi^{u,v}([u]_\Lambda\cap B)}d\mu
=
\int_{[v]_\Lambda\cap B}d\mu\nonumber\\
&=&
\int_{[u]_\Lambda\cap B}{\rm e}^{\psi_f(x,\phi^{u,v}(x))}\,d\mu(x)\nonumber\\
&=&
\int_{[u]_\Lambda\cap B}\exp \sum_{g\in G}\big(f(g\cdot \phi^{u,v}(x))-f(g\cdot x)\big)\,d\mu(x)\,.\nonumber
\end{eqnarray}
Summing over $[v]_\Lambda\in\cW$, and observing that
$\sum_{[v]_\Lambda\in\cW}1_{[v]_\Lambda}(x)=1$ if $x\in A_\cW$,
\begin{eqnarray}\label{eqB}
\sum_{[v]_\Lambda\in\cW}\int_{[v]_\Lambda\cap B}d\mu
&=&
\int 1_B(x)(\sum_{[v]_\Lambda\in\cW}1_{[v]_\Lambda}(x))\,d\mu(x)
\nonumber\\
&=&
\int_{B}d\mu=
\int_{[u]_\Lambda\cap B}\big(\pi^f_\Lambda(x|x)\big)^{-1}\,d\mu(x)\,,
\end{eqnarray}
since for $x\in[u]_\Lambda\cap B$, see \eqref{eqconvenient},
$$
\sum_{[v]_\Lambda\in\cW}\exp \sum_{g\in G}\big(f(g\cdot \phi^{u,v}(x))-f(g\cdot x)\big)
=
\big(\pi^f_\Lambda(x|x)\big)^{-1} \,.
$$

\begin{lem}\label{lemM5.1}
Let $[u]_\Lambda\in \cW$. Given $\varepsilon>0$, there exists $\Lambda^{\prime\prime}\in\cP_f(G)$ such that for all 
$x,y\in [u]_\Lambda\cap A_\cW$, 
$$
d_{\Lambda^{\prime\prime}}(x,y)\leq\varepsilon^*\implies {\rm e}^{-\varepsilon}\leq\Big|
\frac{\pi_\Lambda^f(x|x)}{\pi_\Lambda^f(y|y)}\Big|\leq{\rm e}^\varepsilon\,.
$$
\end{lem}

\prf
Let $x,y\in [u]_\Lambda\cap A_{\cW}$. 
$$
W_\Lambda^x=\{x_w: x_w=\phi_\Lambda^{u,w}(x)\,,\;[w]_\Lambda\in\cW\}
$$
and
$$
W_\Lambda^y=\{y_w: y_w=\phi_\Lambda^{u,w}(y)\,,\;[w]_\Lambda\in\cW\}\,.
$$
Since $x_w\in W_\Lambda^x$ and $y_w\in W_\Lambda^y$,
$$
\sum_{g\in\Lambda^c}|f(g\cdot x)-f(g\cdot x_w)|<\infty\quad\text{and}\quad
\sum_{g\in\Lambda^c}|f(g\cdot y)-f(g\cdot y_w)|<\infty\,.
$$
Given $\varepsilon>0$, there exists $\Lambda^\prime\supset\Lambda$ such that for all 
$[w]_\Lambda\in\cW$,
$$
\sum_{g\not \in\Lambda^\prime}|f(g\cdot x)-f(g\cdot x_w)|\leq\varepsilon/4\;,\quad
\sum_{g\not \in\Lambda^\prime}|f(g\cdot y)-f(g\cdot y_w)|\leq\varepsilon/4\,,
$$
and
\begin{eqnarray*}
& &\Big|\sum_{g\in G}\big(f(g\cdot x_w)-f(g\cdot x)\big)-
\sum_{g\in G}\big(f(g\cdot y_w)-f(g\cdot y)\big)\big)\Big|\leq\\
& &
 \sum_{g\in\Lambda^\prime}\big| f(g\cdot x)-f(g\cdot y)\big|+
 \sum_{g\in\Lambda^\prime}\big| f(g\cdot x_w)-f(g\cdot y_w)\big|+\varepsilon/2\,.
\end{eqnarray*}
Let $\varepsilon^\prime=|\Lambda^\prime|^{-1}\varepsilon/4$ and 
$\delta>0$ such that
$$
d(x,y)\leq\delta\implies |f(g\cdot x)-f(g\cdot y)|\leq\varepsilon^\prime\,\;\forall g\in\Lambda^\prime\,.
$$
For any $\Lambda^{\prime\prime}\supset\Lambda$, if $d_{\Lambda^{\prime\prime}}(x,y)\leq\varepsilon^*$, then 
$d_{\Lambda^{\prime\prime}}(x_w,y_w)\leq\varepsilon^*$.
By lemma \ref{lemdistance} ii) there exists $\Lambda^{\prime\prime}$ such that $d_{\Lambda^{\prime\prime}}(z,z^\prime)\leq\varepsilon^*
\implies d(z,z^\prime)\leq\delta$. Let $d_{\Lambda^{\prime\prime}}(x,y)\leq\varepsilon^*$; uniform continuity of $f$ implies
$$
 \sum_{g\in\Lambda^\prime}\big| f(g\cdot x)-f(g\cdot y)\big|\leq
\varepsilon/4
 \quad\text{and}\quad
  \sum_{g\in\Lambda^\prime}\big| f(g\cdot x_w)-f(g\cdot y_w)\big| \leq 
\varepsilon/4\,.
$$
Since $x,y\in A_\cW$, there exists a bijection between $\cW_\Lambda(x)$ and $\cW_\Lambda(y)$, $x_w\mapsto y_w$.
The lemma follows from \eqref{eqconvenient} and these estimates.
\qed

Let $\{\Lambda_m\}\subset\cP_f(G)$ be an increasing family such that $\bigcup_mA_m=G$ and 
$\Lambda_m\supset\Lambda^{\prime\prime}$ of lemma  \ref{lemM5.1}.  Let $y\in A_\cW$. 
The set $B=[y]_{\Lambda_m\backslash\Lambda}\cap A_\cW$ is measurable and $\Lambda$-saturated. From  \eqref{eqB}
and denoting by $y_u$ the unique element of $[u]_\Lambda\cap W_\Lambda^y$,
\begin{eqnarray*}
\mu(B)=
\int_{[u]_\Lambda\cap B}\big(\pi^f_\Lambda(x|x)\big)^{-1}\,d\mu(x)
&\leq& 
{\rm e}^{\varepsilon}\big(\pi^f_\Lambda(y_u|y_u)\big)^{-1}\mu([u]_\Lambda\cap B)\\
&=&
{\rm e}^{\varepsilon}\big(\pi^f_\Lambda(y_u|y)\big)^{-1}\mu([u]_\Lambda\cap B)
\end{eqnarray*}
since  by \eqref{eqequality} $\pi^f_\Lambda(y_u|y_u)=\pi^f_\Lambda(y_u|y)$.
Similarly
$$
\mu( B)\geq 
{\rm e}^{-\varepsilon}
\big(\pi^f_\Lambda(y_u|y)\big)^{-1}\mu([u]_\Lambda\cap B)\,.
$$
These inequalities can be rewritten as (provided that $\mu(B)>0$)
$$
{\rm e}^{-\varepsilon}\mu([u]_\Lambda| B)\leq 
\pi^f_\Lambda(y_u|y)\leq 
{\rm e}^{\varepsilon}\mu([u]_\Lambda|B)\,.
$$
Let $\varphi\in C(X)$. Define
\begin{eqnarray*}
\varphi^+_m(x)
&:=&
\sup\{\varphi(x^\prime): d_{\Lambda_m}(x^\prime, x)\leq\varepsilon^*\}\\
\varphi^-_m(x)
&:=&
\inf\{\varphi(x^\prime): d_{\Lambda_m}(x^\prime, x)\leq\varepsilon^*\}\,.
\end{eqnarray*}
Then
\begin{eqnarray*}
\frac{{\rm e}^{-\varepsilon}}{\mu(B)}\int_{[u]_\Lambda\cap B}\varphi^-_m(x^\prime)\,d\mu(x^\prime)
&\leq &
\varphi(y_u)\,{\rm e}^{-\varepsilon}\frac{\mu([u]_\Lambda\cap B)}{\mu(B)}\\
&\leq &
\varphi(y_u)\,\pi^f_\Lambda(y_u|y)\\
&\leq &
\varphi(y_u)\,{\rm e}^{\varepsilon}\frac{\mu([u]_\Lambda\cap B)}{\mu(B)}\\
&\leq &
\frac{{\rm e}^{\varepsilon}}{\mu(B)}\int_{[u]_\Lambda\cap B}\varphi^+_m(x^\prime)\,d\mu(x^\prime)\,.
\end{eqnarray*}
Summing over $u\in\cW$ and then integrating  over $B$ with respect to $y$,
$$
{\rm e}^{-\varepsilon}\int_{ B}\varphi_m^-(x^\prime)\,d\mu(x^\prime)\leq 
\int_B\pi^f_\Lambda(\varphi|y)\,d\mu(y)\leq
{\rm e}^{\varepsilon}\int_{ B}\varphi_m^-(x^\prime)\,d\mu(x^\prime)\,.
$$
The partition $\cP$ of $X$ defined by the sets $A_\cW$ is finite. One may assume that the above estimates are true for any $A_\cW$ with the same $m$.
The partition  $\cP_m=\{[y]_{\Lambda_m\backslash\Lambda}\}\vee\{A_\cW\}$ of $X$ is finite. Summing over $B\in\cP_m$ and then letting $m\ra\infty$ one gets
$$
\int_{X}\pi^f_\Lambda(\varphi|x)\,d\mu(x)=\int_X\varphi(x)\,d\mu(x)
$$
since $\varepsilon$ is arbitrary.

\noindent
b) Let $\mu$ be a Gibbs measure for $f$  and $\phi_\Lambda^{u,v}$  be given. 

\begin{lem}\label{lemb}
Let $\Lambda\subset\Lambda^\prime\in\cP_f(G)$. For any $\bar x$,
$\phi_\Lambda^{u,v}$ maps bijectively $W_{\Lambda^\prime}^{\bar x}$ onto $W_{\Lambda^\prime}^{\bar x}$.
\end{lem}
\prf
Suppose that $x\in W_{\Lambda^\prime}^{\bar x}$. Either  $\phi_\Lambda^{u,v}(x)=x$, or $x\in[u]_\Lambda$ and there exist $v\not=u$ and a unique
$y\in [v]_\Lambda\cap W_\Lambda^x$ such that
$\phi_\Lambda^{u,v}(x)=y$. Moreover
$$
d_{\Lambda^c}(x,y)\leq\varepsilon^*\;\text{and}\;d_{{\Lambda^\prime}^c}(x,\bar x)\leq\varepsilon^*
\implies d_{{\Lambda^\prime}^c}(y,\bar x)\leq\varepsilon^*\quad\text{i.e. $y\in W_{\Lambda^\prime}^{\bar x}$}\,.
$$
Interchanging the role of $u$ and $v$, $\phi_\Lambda^{u,v}(y)=x$.
\qed

It follows from lemma \ref{lemb} that for all $w\in W_{\Lambda^\prime}^{\bar x}$
\begin{eqnarray*}
\pi_{\Lambda^\prime}^f(1_{\phi_\Lambda^{u,v}([w]_{\Lambda^\prime})}|\bar x)
&=&\pi_{\Lambda^\prime}^f(\phi_\Lambda^{u,v}(w)|\bar x)\\
&=&
{\rm e}^{\psi_f(w,\phi_\Lambda^{u,v}(w))}\pi_{\Lambda^\prime}^f(w|\bar x)\\
&=&
\pi_{\Lambda^\prime}^f({\rm e}^{\psi_f(\cdot,\phi_\Lambda^{u,v}(\cdot))}1_{[w]_{\Lambda^\prime}}|\bar x)\,.
\end{eqnarray*}
Therefore
\begin{eqnarray*}
\mu\circ\phi_\Lambda^{u,v}([w]_{\Lambda^\prime})
&=&
\int_X\pi_{\Lambda^\prime}^f(1_{\phi_\Lambda^{u,v}([w]_{\Lambda^\prime})}|x)\,d\mu(x)\\
&=&\int_X\pi_{\Lambda^\prime}^f({\rm e}^{\psi_f(\cdot,\phi_\Lambda^{u,v}(\cdot))}1_{[w]_{\Lambda^\prime}}|x)\,d\mu(x)\\
&=&
\int_{[w]_{\Lambda^\prime}}{\rm e}^{\psi_f(x,\phi_\Lambda^{u,v}(x)}\,d\mu(x)\,.
\end{eqnarray*}
This identity holds for any cylinder set $[w]_{\Lambda^\prime}$, $\Lambda^\prime\supset\Lambda$, but also for 
$\Lambda^\prime\subset\Lambda$. Indeed, if $\Lambda^\prime\subset\Lambda$, then one decomposes $[w]_{\Lambda^\prime}$ as
a disjoint union of cylinder sets $[z]_{\Lambda}$ and $\phi_{\Lambda}^{u,v}([w]_{\Lambda^\prime})$ is the disjoint union of the cylinders
$\phi_{\Lambda}^{u,v}([z]_{\Lambda})$ because $\phi_{\Lambda}^{u,v}$ is a bijection.
 Consequently
$$
\frac{d\mu\circ\phi_\Lambda^{u,v}}{d\mu}(x)={\rm e}^{\psi_f(x,\phi_\Lambda^{u,v}(x))}\,.
$$
\qed

\subsection{Dependency on the choice of the expansive constant}\label{independence}

For the dynamical case $(G,X)$ the basic sets $[u]_\Lambda$ and $W_\Lambda^x$ defined in \eqref{WLambda} depend on the choice of $\varepsilon^*$ and consequently 
also the map $\phi_\Lambda^{u,v}$, $M_\Lambda$ and the notions 
of $\Lambda$-saturated sets,  Gibbs conformal  and Gibbs measures depend on the choice of $\varepsilon^*$ . In this subsection one writes this dependency explicitly, e.g. by writing $[u]_\Lambda^\varepsilon$ and $W_\Lambda^{x,\varepsilon}$ instead of $[u]_\Lambda$ and $W_\Lambda^x$ if $\varepsilon^*=\varepsilon$.

\begin{pro}\label{proM}
Let $0<\varepsilon^\prime<\varepsilon$ be  expansive constants. The measure
$\mu$ is a Gibbs conformal measure for $\varepsilon^*=\varepsilon$ if and only if it is a Gibbs conformal measure for $\varepsilon^*=\varepsilon^\prime$.
\end{pro}

\prf 
a) Suppose that $\mu$ is Gibbs conformal for $\varepsilon^*=\varepsilon$.
Let $\varepsilon^\prime<\varepsilon$, $x\in [u]^{\varepsilon^\prime}_\Lambda$ and suppose that 
$[v]^{\varepsilon^\prime}_\Lambda\cap W_\Lambda^{x,\varepsilon^\prime}=\{y\}$. Then
$[v]^{\varepsilon}_\Lambda\cap W_\Lambda^{x,\varepsilon}=\{y\}$ and
$$
\phi_\Lambda^{u,v,\varepsilon^\prime}(x)=\phi_\Lambda^{u,v,\varepsilon}(x)=y\,.
$$
Let $A\subset X$ and decompose $A$ into $A_1$ and $A_2$,
$$
A_1:=\{x\in A\colon \phi_\Lambda^{u,v,\varepsilon^\prime}(x)=x\}\quad\text{and}\quad
A_2:=\{x\in A\colon \phi_\Lambda^{u,v,\varepsilon^\prime}(x)\not=x\}\,.
$$
Then 
$$
\phi_\Lambda^{u,v,\varepsilon^\prime}(A)=A_1\sqcup \phi_\Lambda^{u,v,\varepsilon^\prime}(A_2)=
A_1\sqcup \phi_\Lambda^{u,v,\varepsilon}(A_2)
$$
and
\begin{eqnarray*}
\mu\circ\phi_\Lambda^{u,v,\varepsilon^\prime}(A)
&=&
\mu(A_1)+\mu(\phi_\Lambda^{u,v,\varepsilon}(A_2))\\
&=&
\int_{A_1}\,d\mu+\int_{A_2}\exp\psi_f(x,\phi_\Lambda^{u,v,\varepsilon}(x))\,d\mu(x)\\
&=&
\int_{A}\exp\psi_f(x,\phi_\Lambda^{u,v,\varepsilon^\prime}(x))\,d\mu(x)\,,
\end{eqnarray*}
since $\psi_f(x,\phi_\Lambda^{u,v,\varepsilon^\prime}(x))=0$ on $A_1$.
Therefore $\mu$ is a Gibbs conformal measure for $\varepsilon^*=\varepsilon^\prime$.\\
b) Let $\mu$ be a  Gibbs conformal measure for $\varepsilon^*=\varepsilon^\prime$. 
Consider  $\phi_\Lambda^{u,v,\varepsilon}$ with
 $[u]_\Lambda^\varepsilon\cap [v]_\Lambda^\varepsilon=\emptyset$.
Let $\bar\Lambda\supset\Lambda$ such that $d_{\Lambda^c}(x,y)\leq\varepsilon^*$ implies $d_{\bar\Lambda^c}(x,y)\leq\varepsilon^\prime$
(lemma \ref{lemdistance}). In particular $W_{\bar\Lambda}^{x,\varepsilon^\prime}\supset W_{\Lambda}^{x,\varepsilon}$.
The subsets $[u]_\Lambda^\varepsilon$, respectively $[v]_\Lambda^\varepsilon$, can be partitioned into subsets
$[u_i]_{\bar\Lambda}^{\varepsilon^\prime}$, respectively $[v_j]_{\bar\Lambda}^{\varepsilon^\prime}$. The set $[u]_\Lambda^\varepsilon$ is decomposed into
$$
A_1=\{x\in [u]_\Lambda^\varepsilon: \phi_\Lambda^{u,v,\varepsilon}(x)=x\}\quad\text{and}\quad
A_2=\{x\in [u]_\Lambda^\varepsilon: \phi_\Lambda^{u,v,\varepsilon}(x)=y\,,\;y\in [v]_\Lambda^\varepsilon\}\,.
$$
Let $x\in A_2\cap[u_i]_{\bar\Lambda}^{\varepsilon^\prime}$. Since $W_{\bar\Lambda}^{x,\varepsilon^\prime}\supset W_{\Lambda}^{x,\varepsilon}$,
$$
y\in  [v_j]_{\bar\Lambda}^{\varepsilon^\prime}\cap W_{\bar\Lambda}^{x,\varepsilon^\prime}\quad\text{for some $j$ and }
\phi_\Lambda^{u,v,\varepsilon}(x)=\phi_{\bar\Lambda}^{u_i,v_j,\varepsilon^\prime}(x)\,.
$$
Let
$$
A_2^{(i,j)}:=
\{x\in A_2\cap[u_i]_{\bar\Lambda}^{\varepsilon^\prime}\colon 
\phi_{\bar\Lambda}^{u_i,v_j,\varepsilon^\prime}(x)\in [v_j]_{\bar\Lambda}^{\varepsilon^\prime}\}\,.
$$
Therefore,
$$
\phi_\Lambda^{u,v,\varepsilon}(A_2)
=\bigsqcup_{(i,j)}\phi_{\bar\Lambda}^{u_i,v_j,\varepsilon^\prime}(A_2^{(i,j)})
$$
and
\begin{eqnarray*}
\mu(\phi_{\Lambda}^{u,v,\varepsilon}([u]_\Lambda^\varepsilon))&=&\mu(A_1)+\mu(\phi_{\Lambda}^{u,v,\varepsilon}(A_2))\\
&=&
\mu(A_1)+\sum_{(i,j)}
\mu(\phi_{\Lambda}^{u_i,v_j,\varepsilon^\prime}(A_2^{(i,j)}))\\
&=&
\mu(A_1)+\sum_{(i,j)}
\int_{A_2^{(i,j)}}\exp\psi_f(x,\phi_{\Lambda}^{u_i,v_j,\varepsilon^\prime}(x))\,d\mu\\
&=&
\mu(A_1)+\int_{A_2}\exp\psi_f(x,\phi_{\Lambda}^{u,v,\varepsilon}(x))\,d\mu\\
&=&
\int_{[u]_\Lambda^\varepsilon}\exp\psi_f(x,\phi_{\Lambda}^{u,v,\varepsilon}(x))\,d\mu\,.
\end{eqnarray*}
Hence $\mu$ is Gibbs conformal  for $\varepsilon^*=\varepsilon$.
\qed

\section{$\beta$-shifts}\label{sectionbeta}
\setcounter{equation}{0}

Almost all proofs of  existence of  Gibbs measures in the literature are based on the property that $M_\Lambda$ maps $C(X)$ into $C(X)$
(proposition \ref{procontinuityMLambda}).
For  the natural extension of a $\beta$-shift the set of $\beta>1$ such that proposition \ref{procontinuityMLambda} holds is the class $C_3$ in \cite{Sc}, which has Hausdorff dimension
$1$, but which is of Lebesgue measure $0$.  To prove the 
main theorem of this section, theorem \ref{thmbeta}, one relies on theorem \ref{thmM5.2}. 
For some $\beta$ an equilibrium measure is not a weak Gibbs measure \cite{PS}, hence for those $\beta$'s one has examples of Gibbs measures which are not weak Gibbs measures. 

\subsection{Natural extension of $\beta$-shifts}\label{subsectionnaturalextension}
 
Let $\beta>1$ be fixed. The case $\beta\in\N$ is special and corresponds to the full shift.  From now on 
$\beta\not\in\N$.
For $t\in\R$, let
$ \lceil t\rceil :=\min\{i\in\Z\colon i\geq t\}$.
One defines $\tb:=\lceil \beta\rceil$.  Consider the
$\beta$-expansion of $1$,
$$
1=\sum_{i=1}^\infty c(i)\beta^{-i}\,,
$$
which is  given by the algorithm
$$
  r_0:=1,\; c(i+1):=\lceil \beta\, r_i\rceil-1,\;
  r_{i+1}:=\beta \,r_i-c(i+1),\;i\in\Z_+\,,
$$
which insures that $r_i>0$ for all $i\in\Z_+$. It follows that $c(1)=\lceil\beta\rceil-1>0$ and
$c^\beta:=(c(1),c(2),\ldots)$ cannot end with zeros only. For sequences
$(a(1),a(2),\ldots)$ and $(b(1),b(2),\ldots)$ the lexicographical
order is defined by $(a(1),a(2),\ldots)\prec (b(1),b(2),\ldots)$ if and only if for the
least index $i$ with $a(i)\neq b(i)$, $a(i)<b(i)$. Let
$\tA:=\{0,\ldots,\tb-1\}$; the one-sided $\beta$-shift is
$$
  X^\beta:=\left\{ x=(x(1),x(2),\ldots)\colon  x(i)\in\tA,\;
    T^kx \preceq c^\beta\;\forall \,k\in\Z_+\right\}\,,
 $$
where $T$ is the left shift operator.
In particular $T^kc^\beta \preceq c^\beta$ for all $k\in\Z_+$,  so that $X^\beta$ is
a shift-invariant closed subset of $\tA^\N$ (with product topology). The language of the shift $X^\beta$ is denoted by $\cL^\beta$ 
and the set of the words of length $n$ by
$\cL^\beta_n$; the empty-word is $\epsilon$, $\cL^\beta_0=\{\epsilon\}$. The length of a word $w$ is $|w|$. 
A word\index{Word!} 
$w$ of length $n$ is also written  $w(1)\cdots w(n)$ as the concatenation of the $n$ words $w(j)$ of length $1$.

The shift space $X^\beta$ can be described by a labeled graph $\cG^\beta=(\cV,\cE^\beta)$ where $\cV:=\{q_j\colon j\in \Z_+\}$.
The set of vertices of the graph is into $1-1$-correspondence with the set of the prefixes of $c^\beta$, including
$\epsilon$: $q_0$ corresponds to $\epsilon$ and $q_m$, $m\geq 1$, corresponds to the
prefix $c(1)\cdots c(m)$ of $c^\beta$.  
The root of the graph is the vertex $q_0$.
There is an edge $q_0\ra q_0$, labeled by $k$, for each $k=0,\ldots,\tb-2$, and there is an edge $q_{j-1}\ra q_ j$ labeled by $c(j)$  for each $j\in\N$.
If the label $c(j)$ of $q_{j-1}\ra q_j$ is different from $0$, then there are $c(j)$ edges $q_{j-1}\ra q_0$ labeled by $0,\ldots, c(j)-1$.
Each word $w(1)\cdots w(n)\in\cL^\beta$ can always be presented by a path of length $n$ in $\cG^\beta$ starting with vertex $q_0$. For a word $w\in\cL^\beta$ 
one defines $q(w)$ as the end vertex of this path starting at $q_0$ and presenting $w$.
One can concatenate two words  $w$ and $w^\prime$ if and only if there is a path $\eta$ in $\cG^\beta$ presenting $w$ and a path $\eta^\prime$ in 
$\cG^\beta$ presenting $w^\prime$, so that $\eta$ ends at vertex $q$ and $\eta^\prime$ starts at vertex $q$.

 Let $u$ be a prefix of $c^\beta$. Set
\begin{equation}\label{eq5.4.1}
z^\beta(u):=\begin{cases}
p& \text{if $u=c(1)\cdots c(\ell)$, $c(\ell+1)=\cdots =c(\ell+p)=0$, $c(\ell+p+1)>0$}\\
0 &\text{if $u=\epsilon$.}
\end{cases}
\end{equation}
For each prefix $u$ of $c^\beta$ one defines a new word $\widehat{u}$ with $q(\widehat{u})=q_0$ and which differs from $u$ by a single letter.
If $u=c(1)\cdots c(n)$ let $c(j)$ be the last letter in $c(1)\cdots c(n)$ which is different from $0$, and set
\begin{equation}\label{eq5.4.3}
\widehat{c}(\ell):=\begin{cases}
c(\ell) & \text{if $\ell\not=j$}\\
c(j)-1 &\text{if $\ell=j$}
\end{cases}\quad\text{and}\quad\hat{u}:=\widehat{c}(1)\cdots\widehat{c}(n)\,.
\end{equation}
For any word $w\in\cL^\beta$ there is a unique decomposition of $w$ into 
\begin{equation}\label{eqs(v)}
w=vs(w)\quad\text{where $s(w)$ is the largest suffix of $w$, which is a prefix of $c^\beta$.}
\end{equation}
The definition \eqref{eq5.4.1} is extended to any word $w$ by setting
\begin{equation}\label{eq5.4.1bis}
z^\beta(w):=z^\beta(s(w))\,,
\end{equation}
and the transformation $u\mapsto\widehat{u}$ is extended  to any word by setting
\begin{equation}\label{eq5.4.4}
\widehat{w}:=\begin{cases}
w &\text{if $s(w)=\epsilon$}\\
v\widehat{u} &\text{if $s(w)=u$.}
\end{cases}
\end{equation}
By convention $\widehat{\epsilon}:=\epsilon$ 

\begin{lem}\label{lem5.4.1}
i) If $a$ and $b$ are prefixes of $c^\beta$ and $ab\in\cL^\beta$, then $ab$ is a prefix of $c^\beta$.\\
ii) Let $w=vu$, $s(w)=u$. Then $q(v)=q_0$ and $s(\widehat{w})=s(v\widehat{u})=\epsilon$.\\
iii) $q(w)\not=q_0$ if and only if $s(w)\not=\epsilon$.\\
iv) There is a constant $q$ such that  the map on $\cL^\beta$, $w\mapsto\widehat{w}$, is at most
$q$-to-$1$.\\
v) Let $w\in\cL^\beta_n$ and $u$ a prefix of $c^\beta$. If $uw\in\cL^\beta$, then
$z^\beta (uw)\geq z^\beta(w)$.\\
vi) Any prefix $u$ of $c^\beta$ can be extended to the right by $0$.
\end{lem}

\prf
Statements i) to iv) are proved in \cite{PS} (lemma 1 and lemma 2).
Suppose that $w=vs(w)$, so that $z^\beta(w)=z^\beta(s(w))$. 
Since $s(w)$ is maximal, either $s(uw)=s(w)$ or $|s(uw)|>|s(w)|$ and
by maximality of $s(w)$ and statement i) $uw$ is a prefix of $c^\beta$.
In particular $uv$ and $s(w)0\cdots 0$ ($|z^\beta(s(w))|$ letters $0$) are prefixes of $c^\beta$ and again by i)
$uvs(w)0\cdots 0$ is a prefix of $c^\beta$, so that
$z^\beta (u w)\geq z^\beta(w)$.\\
vi) Let $u=c(1)\cdots c(n)$. If $c(n+1)>0$, then one changes $c(n+1)$ into $0$ and then extend this word by $0$ only.
If $c(n+1)=\cdots =c(n+p)=0$ and $c(n+p+1)>0$, then one changes $c(n+p+1)$ into $0$ and then extend this word by $0$ only.
\qed

\begin{defn}
The   the natural extension $\Sigma^\beta$ of $X^\beta$ is
$$
\Sigma^\beta=\{x\in\tA^\Z\colon \forall k\in\Z\,,\;(x(k),x(k+1),\ldots)\in X^\beta\}\,.
$$
The group $\Z$ acts on $\Sigma^\beta$ by the left shit operator $T$. 
\end{defn}

Let $\Lambda=[k,\ell]$ be an interval in $\Z$ and $J_\Lambda$ the projection operator
$$
x\mapsto J_\Lambda(x):=x(k)\cdots x(\ell)\equiv x_\Lambda\,.
$$
Set $\Sigma^\beta_\Lambda:=J_\Lambda(\Sigma^\beta)$.
One decomposes $x$ into $x_-x_\Lambda x_+$, where $x_\Lambda\in\Sigma^\beta_\Lambda$,
$x_-(j)=x(j)$, $j<k$, and $x_+(j)=x(j)$, $j>\ell$. If
$J_\Lambda(x)=w$, then $\overline{w}$ is the element of $\Sigma^\beta$ defined by
$$
J_\Lambda(\overline{w})=w\quad \text{and}\quad \overline{w}(j)=0\,,\;\forall j\not\in [k,\ell]\,.
$$
Existence of $\overline{w}$ is a consequence of lemma \ref{lem5.4.1} vi. For $a\in \Sigma_\Lambda^\beta$,
$$
[a]:=\{x\in\Sigma^\beta\colon J_\Lambda(x)=a\}\,.
$$
On $\Sigma^\beta$ one chooses the distance defined by
\begin{equation}\label{distance}
d(x,y)=\epsilon_m\quad\text{if}\;m=\min\{n\colon x(k)\not =y(k)\,,\;\text{$k=n$ or $k=-n$}\}\,,
\end{equation}
where  $\epsilon_n$, $n\geq 0$, is a strictly decreasing sequence, 
$\epsilon_0=1$ and $\lim_n\epsilon_n=0$.

\subsection{An equilibrium measure is a Gibbs measure}\label{subsectionBowen}

Let $f\in C(\Sigma^\beta)$ and $\Lambda\subset\Z$ be an interval.
$$
\Xi_\Lambda(f):=\sum_{w\in \Sigma_{\Lambda}^\beta}\exp\sum_{j\in \Lambda}f(T^j\overline{w})
\quad\text{and}\quad P_\Lambda(f):=\frac{1}{|\Lambda|}\ln\Xi_\Lambda(f)\,.
$$
The pressure $P(f)$ is equal to 
\begin{equation}\label{pressure}
P(f)=\lim_{n\ra\infty}P_{[-n,n]}(f)\,.
\end{equation}

\begin{defn}\label{defnBowen}
A  continuous function $f$ verifies Bowen's condition if there exists a constant $C$ such that for all intervals $\Lambda\subset\Z$
$$
\sup\big\{\big|\sum_{k\in\Lambda}(f(T^kx)-f(T^ky))|: J_\Lambda(x)=J_\Lambda(y)\big\}\leq C\,.
$$
\end{defn}

\noindent
{\bf Remark.\,}  A continuous function has \emph{summable variation} if 
$$
{\rm var}_n(f):=\sup\big\{|f(x)-f(y)|\colon x(k)=y(k)\;\forall k\in [0,n-1]\big\}
$$
and $\sum_{n\geq 1} {\rm var}_n(f)<\infty$.
Such a function verifies Bowen's condition, lemma 1.15  \cite{Bo}. They are also  in  $\B(\Sigma^\beta)$.
\qed

\begin{thm}[\cite{CT}]\label{thmCT}
For any $\beta >1$, if $f\in C(\Sigma^\beta)$ verifies Bowen's condition, then there is a unique equilibrium measure
$\nu$  for $f$.  Moreover
$$
\limsup_n\frac{1}{n}\sum_{k=0}^{n-1}f(T^kc^\beta)<P(f)\,.
$$
\end{thm}

\prf  
There is a $1$-to-$1$ map from the $T$-invariant measures on $\Sigma^\beta$ and the $T$-invariant measures on
$X^\beta$, which preserves the entropy and the pressure. For $X^\beta$ theorem \ref{thmCT} is a consequence of proposition 3.1 and theorem C
in \cite{CT}.
\qed

\begin{thm}\label{thmbeta}
Suppose that $f\in \B(\Sigma^\beta)$ verifies Bowen's condition.
Then for any $\beta >1$ the unique equilibrium measure $\nu$ for $f$ is a Gibbs measure.
\end{thm}

Theorem \ref{thmbeta} follows from theorem  \ref{thmM5.2} if one proves that the equilibrium measure $\nu$ is Gibbs conformal.

\begin{pro}\label{pronuf}
Suppose that $f\in C(\Sigma^\beta)$ verifies Bowen's condition and that $\Lambda$ is an interval. Set
\begin{equation}\label{mun}
\mu_\Lambda:=\frac{\sum_{w\in \Sigma^\beta_{\Lambda}}\exp\big(\sum_{i\in \Lambda}f(T^i\overline{w})\big)\delta_{\overline{w}}}{\Xi_\Lambda(f)}\,.
\end{equation}
Then  the unique equilibrium measure $\nu=\lim_n\nu_n$ where
$$
\nu_n:=\frac{1}{2n+1}\sum_{j=-n}^nT^{-j}\mu_{[-n,n]}\,.
$$
\end{pro}

\prf
A $T$-invariant probability measure $\nu$ is an equilibrium measure for a continuous function $f$ if and only if 
$\nu$ is a tangent functional to the pressure $P$ at $f$ (see \cite{Wa}  theorems 8.2 and 9.5), 
$$
P(f+\varphi)\geq P(f)+\int \varphi\,d\nu\,,\;\text{for all continuous functions $\varphi$}\,.
$$
When there is a unique tangent functional $\nu$ to the pressure at $f$ one can express
$\nu([u])$, $u\in\cL^\beta$,  using a classical result about
differentiability of a convex function, here the pressure, which is a pointwise limit of convex functions, theorem 25.7 in \cite{Ro}.
Let $u\in\cL^\beta$, $|u|=k$, and set
\begin{equation*}
I_{u}(y):=\left\{\begin{array}{lll}
1 &\text{if $y_0\cdots y_{k-1}= u$}\\
0 &\text{otherwise.}
\end{array}\right.
\end{equation*}
Then
\begin{eqnarray}\label{eq13}
  \nu(I_{u})
  &=&
  \left .\frac {d}{dt} \lim_{n\to\infty}
  P_{[-n,n]}(f+t\,I_{u})\right|_{t=0}=
 \lim_{n\to\infty}\left .\frac {d}{dt} P_{[-n,n]}(f+t\,I_{u})\right|_{t=0}\nonumber\\
 &=&
\lim_{n\ra\infty} \frac{1}{2n+1}\sum_{j=-n}^n
\frac{\sum_{w\in \Sigma_{[-n,n]}^\beta}I_{u}(T^j \overline{w}) \exp\sum_{i\in[-n,n]} f(T^i\overline{w})}
{ \Xi_{[-n,n]}(f)}\nonumber\\
&=&
\lim_{n\ra\infty} \frac{1}{2n+1}\sum_{j=-n}^n\mu_{[-n,n]}(I_u\circ T^j)\,.\nonumber
\end{eqnarray}
\qed

\subsection{Proof of theorem \ref{thmbeta}}

For convenience one writes  $\Sigma$ for $\Sigma^\beta$, $\Sigma_\Lambda$ for $\Sigma^\beta_\Lambda$, $\cL$ for $\cL^\beta$.
Let $\Lambda\subset\Z$ be an interval and $x\in\Sigma$. If $J_\Lambda(x)=w$, then on writes $x$ as $x_-wx_+$.
Thus for $u,v\in \Sigma_\Lambda$,
\begin{equation}\label{eqinvolutionbeta}
\phi^{u,v}_\Lambda(y):=\begin{cases}
x_-vx_+&\text{if $y=x_-ux_+\in\Sigma$ and $x_-vx_+\in\Sigma$}\\
x_-ux_+ &\text{if $y=x_-vx_+\in\Sigma$ and $x_-ux_+\in\Sigma$}\\
y &\text{otherwise.}
\end{cases}
\end{equation}
One defines a tree $\cT(u)$ describing all $x\in [u]$. The root of the tree is $u$. 
The vertices of level $k$ of the tree are the pairs $(w_-,w_+)$, $|w_-|=|w_+|=k$, such that
there exists $x\in[u]$, $x=\cdots w_-uw_+\cdots$.
In other words, $x\in[u]$ if and only if  $x_-ux_+$ corresponds to an infinite path of the tree. 
If for all $(w_-,w_+)$ along a path of the tree $\cT(u)$ the words $w_-vw_+\in\cL$, then $x_-vx_+$ corresponds to an infinite path of $\cT(b)$, 
and $x_-vx_+\in [v]$.
By definition 
$$
\cT(u,v):=\{(w_-,w_+)\colon \text{$(w_-,w_+)$ is a common vertex of $\cT(u)$ and  $\cT(v)$}\}\,.
$$
Thus
$$
\phi^{u,v}_\Lambda(x_-ux_+)=x_-vx_+
$$
if and only if 
$x_-ux_+$ and $x_-vx_+$ correspond to paths of $\cT(u)$ and $\cT(v)$. On the other hand
$$
\phi^{u,v}_\Lambda(x_-ux_+)=x_-ux_+
$$
if and only if   $x_-ux_+\in\Sigma$ and $x_-vx_+\not\in\Sigma$, that is, if and only if there exists a vertex $(w_-,w_+)$ in the tree $\cT(u)$ so that 
$w_-vw_+\not\in\cL$. 
Let 
$$
V(u):=\{\text{$(w_-,w_+)$ is a vertex of $\cT(u)$ such that $w_-vw_+\not\in \cL$}\}\,.
$$ 
For $x\in[u]$,
\begin{equation}\label{eqouvert}
\phi^{u,v}_\Lambda(x)=x \iff x\in \bigcup_{(w_-,w_+)\in V(u)}[w_-uw_+]\,.
\end{equation}
Consequently $[u]$ is decomposed into an open subset $A_1$ and a closed subset $A_2$,
$$
A_1:=\{x\in[u]:\phi^{u,v}_\Lambda(x)=x\}\quad\text{and}\quad  A_2:=\{x\in[u]:\phi^{u,v}_\Lambda(x)=x_-vx_+\}\,.
$$
One can decompose $[v]$ in a similar way into
$B_1\sqcup B_2$;  there is a bijection between $A_2$ and $B_2$, and
$$
\phi^{u,v}_\Lambda([u])=A_1\sqcup B_2\,.
$$
Set
\begin{equation}\label{eqA2n}
A_2(n):=\bigcup_{\substack{(w_-,w_+)\in \cT(u,v):\\ |w_-|=|w_+|=n}}[w_-uw_+]\quad\text{and}\quad 
B_2(n):=\bigcup_{\substack{(w_-,w_+)\in \cT(u,v):\\ |w_-|=|w_+|=n}}[w_-vw_+]\,.
\end{equation}
The sequences of sets $A_2(n)$ and $B_2(n)$ are decreasing and 
$$
A_2=\bigcap_nA_2(n)\quad\text{and}\quad B_2=\bigcap_nB_2(n)\,.
$$
Therefore $\nu(A_2)=\lim_n\nu(A_2(n))$. 
Since the set $A_2(n)$ is a  finite disjoint union of cylinder sets,
$$
\nu(A_2(n))=\sum_{\substack{(w_-,w_+)\in \cT(u,v):\\ |w_-|=|w_+|=n}}\nu([w_-uw_+])=
\lim_{m\ra\infty}\sum_{\substack{(w_-,w_+)\in \cT(u,v):\\ |w_-|=|w_+|=n}}\nu_m([w_-uw_+])\,.\nonumber
$$
Similarly,
$$
\nu(B_2(n))=
\sum_{\substack{(w_-,w_+)\in \cT(a,b):\\ |w_-|=|w_+|=n}}\nu([w_-vw_+])
=
\lim_{m\ra\infty}\sum_{\substack{(w_-,w_+)\in \cT(a,b):\\ |w_-|=|w_+|=n}}\nu_m([w_-vw_+])\,.
$$

\begin{lem}\label{lemsuficient}
Let $x\in A_2(n)$ and $\Lambda=[k,\ell]$.
Assume  that $[-m,m]\supset [k-n,\ell+n]$ and $J_{[-m,m]}(x)=aw_-uw_+b$ with $(w_-,w_+)$ a vertex of $\cT(u,v)$,
$|w_-|=|w_+|=n$.
Sufficient conditions so that $aw_-vw_+b\in\cL$ are 
$$
|s(aw_-)|\leq n\quad\text{and}\quad |s(aw_-vw_+)|\leq n\,.
$$
\end{lem}

\prf
By hypothesis $w_-vw_+\in\cL$. The word $aw_-$ is decomposed as in \eqref{eqs(v)}, 
$$
aw_-=rs(aw_-)\,.
$$
If $|s(aw_-)|\leq n$, then $|a|\leq |r|$ and $s(aw_-)$ is a factor of $w_-$.
Since $w_-vw_+\in\cL$, also $aw_-vw_+\in\cL$. 
Similarly, the word $aw_-vw_+$ is decomposed into $ts(aw_-vw_+)$.
If $|s(aw_-vw_+)|\leq n$, then $s(aw_-vw_+)$ is a factor of $w_+$. Since $w_+b$ is a factor of 
$aw_-uw_+b\in\cL$, the word $aw_-vw_+b\in\cL$.
\qed

\begin{lem}\label{lemnotb}
Under the assumptions of lemma \ref{lemsuficient}, let
\begin{eqnarray*}
A_2(n|m)&:=&\{x\in A_2(n): J_{[-m,m]}(x)=aw_-uw_+b\,\;\text{and}\,\;   aw_-vw_+b\not\in\cL\}\\
B_2(n|m)&:=&\{x\in B_2(n): J_{[-m,m]}(x)=aw_-vw_+b\,\;\text{and}\,\;   aw_-uw_+b\not\in\cL\}\,.
\end{eqnarray*}
Then there exist constants $C^\prime<\infty$ and $\kappa>0$ such that (for $m$ and $n$ large enough)
$$
\mu_{[-m,m]}(A_2(n|m))\leq C^\prime{\rm e}^{-\kappa n}\quad\text{and}\quad \mu_{[-m,m]}(B_2(n|m))\leq C^\prime{\rm e}^{-\kappa n}\,.
$$
\end{lem}

\prf
For $x\in A_2(n)$ let $J_{[-m,m]}(x)= aw_-uw_+b$.
If $aw_-vw_+b\not\in\cL$,  then $|s(aw_-)|>n$ or $aw_-vw_+\in\cL$ and  $|s(aw_-vw_+)|> n$.
Therefore $A_2(n|m)\subset A_2^-(n|m)\cup A_2^+(n|m)$,
$$
A_2^-(n|m):=\{x\in A_2(n): |s(aw_-)|>n\}\,,
$$
and
$$
A_2^+(n|m):=\{x\in A_2(n): aw_-vw_+\in\cL\,,\;|s(aw_-vw_+)|>n\}\,.
$$
Suppose that the interval $[-m,m]$ is decomposed into three intervals, 
$$
[-m,m]=\Lambda_-\cup\Lambda_0\cup\Lambda_+\,.
$$ 
Write $w\in\Sigma_{[-m,m]}$ as $w=w_-w_0w_+$, $w_-\in\Sigma_{\Lambda_-}$,
$w_0\in\Sigma_{\Lambda_0}$ and $w_+\in\Sigma_{\Lambda_+}$.
Then
\begin{eqnarray*}
\sum_{j\in[-m,n]}f(T^j\overline{w})
&=&
\sum_{j\in\Lambda_-}\big(f(T^j \overline{w})-f(T^j\overline{w_-})\big)+\sum_{j\in\Lambda_-}f(T^j\overline{w_-})\\
&+&
\sum_{j\in\Lambda_0}\big(f(T^j \overline{w})-f(T^j\overline{w_0})\big)+\sum_{j\in\Lambda_0}f(T^j\overline{w_0})\\
&+&
\sum_{j\in\Lambda_+}\big(f(T^j \overline{w})-f(T^j\overline{w_+})\big)+\sum_{j\in\Lambda_+}f(T^j\overline{w_+})\,.
\end{eqnarray*}
Since $f$ verifies Bowen's condition,
$$
\Big|\sum_{j\in[-m,n]}f(T^j\overline{w})-\Big(\sum_{j\in\Lambda_-}f(T^j\overline{w_-})
+\sum_{j\in\Lambda_0}f(T^j\overline{w_0})+\sum_{j\in\Lambda_+}f(T^j\overline{w_+})\Big)\Big|\leq 3C\,.
$$
For an interval $\Lambda$ set
$$
\Xi_\Lambda^0(f):=\sum_{\substack{w\in\Sigma_\Lambda:\\q(w)=q_0}}\exp\sum_{j\in\Lambda}f(T^j\overline{w})\,.
$$
Since the map $w\mapsto\hat{w}$ is at most $q$-to-one,
$$
\Xi_\Lambda^0(f)\leq \Xi_\Lambda(f)\leq q{\rm e}^{2C+2\|f\|_\infty}\,\Xi_\Lambda^0(f)\,.
$$
Therefore, if $[-m,m]=\Lambda_-\cup\Lambda_0\cup\Lambda_+$,
\begin{eqnarray}\label{mumw0}
\mu_{[-m,m]}([w_0])
&\leq&
{\rm e}^{6C}
\frac{\Xi_{\Lambda_-}(f)\Xi_{\Lambda_+}(f)\Xi_{\Lambda_0}(f)}
{\Xi_{\Lambda_-}^0(f)\Xi_{\Lambda_0}^0(f)\Xi_{\Lambda_+}^0(f)}
\cdot\frac{\exp\sum_{k\in\Lambda_0}f(T^k\overline{w_0})}{\Xi_{\Lambda_0}(f)}\nonumber\\
&\leq&
q^3{\rm e}^{12C+6\|f\|_\infty}\cdot\mu_{\Lambda_0}([w_0])\,.
\end{eqnarray}
From \eqref{mumw0} and  theorem \ref{thmCT} one gets the existence of $\kappa>0$ such that for $|\Lambda_0|$ large enough, 
\begin{equation}\label{fund}
\text{$w_0$ prefix of $c^\beta$} \implies \mu_{\Lambda_0}([w_0])\leq{\rm e}^{-\kappa |\Lambda_0|}\,.
\end{equation}
Let $c_j$ be the prefix of $c^\beta$ of length $j$. One decomposes $[-m,m]$ into
$$
[-m, k-j)\cup[k-j,k)\cup[k,m]\,.
$$
If $w=a w_-uw_+b$ and $s(aw_-)=c_j$,  then $J_{[k-j,k)}(\bar w)=c_j$ and ($m$ and $n$ large enough)
$$
\mu_{[-m,m]}(A_2^-(n|m))\leq\sum_{j>n}q^3{\rm e}^{12C+6\|f\|_\infty}{\rm e}^{-\kappa j}\,.
$$
A similar proof gives
$$
\mu_m(A_2^+(n|m))\leq\sum_{j>n}q^3{\rm e}^{12C+6\|f\|_\infty}{\rm e}^{-\kappa j}\,.
$$
\qed

\begin{lem}\label{lemJ0}
Let $[k-n,\ell+n]\supset \Lambda^{\prime\prime}:=[k-t,\ell+t]$ and $x\in A_2(n)$. Then $\overline{J_{\Lambda^{\prime\prime}}(x)}\in A_2$.
\end{lem}

\prf
By hypothesis $x\in [w_-uw_+]$ for some $(w_-,w_+)\in\cT(u,v)$, $|w_-|=|w_+|=n$.
Therefore $w_-uw_+\in\cL$ and $w_-vw_+\in\cL$, so that $J_{\Lambda^{\prime\prime}}(w_-uw_+)\in\cL$ and 
$J_{\Lambda^{\prime\prime}}(w_-vw_+)\in\cL$. Hence  $\overline{J_{\Lambda^{\prime\prime}}(x)}=
\overline{J_{\Lambda^{\prime\prime}}(w_-uw_+)}$ and
$$
\phi_\Lambda^{u,v}(\overline{J_{\Lambda^{\prime\prime}}(w_-uw_+)})=\overline{J_{\Lambda^{\prime\prime}}(w_-vw_+)}\,.
$$
\qed

Let $\varepsilon>0$, $\Lambda=[k,\ell]$  and  $\Lambda^\prime:=[k-r,\ell+r]$. Since $f\in\B(\Sigma)$,  there exists $r$ such that
\begin{equation}\label{eqr}
\sup\big\{\sum_{j\not\in\Lambda^\prime}|f(T^jy)-f(T^jx)|: x(s)=y(s)\,\;\forall s\in\Lambda^c\}\leq \varepsilon\,.
\end{equation}
Let $t$ be large enough, so that $\Lambda^{\prime\prime}=[k-t,\ell+t]\supset\Lambda^\prime$.
Define 
$$
\psi_{f,r}(x):=\sum_{j\in \Lambda^\prime}\big(f(T^j\phi_{\Lambda}^{u,v}(\overline{J_{\Lambda^{\prime\prime}}(x)})-
f(T^j\overline{J_{\Lambda^{\prime\prime}}(x)})\big)\,.
$$

\begin{lem}\label{lemJC}
i) The function $\psi_{f,r}$ is continuous on $\Sigma$. \\
ii) If $x\in A_2$, then
$$
\lim_{t\ra\infty}\psi_{f,r}(x)=\sum_{j\in \Lambda^\prime}\big(f(T^j\phi_{\Lambda}^{u,v}(x))-f(T^jx)\big)\,.
$$
iii) Let $\varepsilon^\prime>0$. There exists $t^\prime$ such that for $n\geq  t\geq t^\prime$ and  $[-m,m]\supset [k-n,\ell+n]$, 
$$
x\in A_2(n)\backslash A_2(n|m)\implies 
\big|\sum_{j\in \Lambda^\prime}\big(f(T^j\phi_{\Lambda}^{u,v}(x))-f(T^jx)\big)-\psi_{f,r}(x)\big|\leq2\varepsilon^\prime\,.
$$
iv) Under the same hypothesis,
$$
\big|\psi_f(x)-\psi_{f,r}(x)\big|\leq \varepsilon+2\varepsilon^\prime\,.
$$
\end{lem}

\prf
i) The continuity of $\psi_{f,r}$ on $\Sigma$ is a consequence of
$$
\psi_{f,r}(x)=\psi_{f,r}(y)\quad\text{if $J_{\Lambda^{\prime\prime}}(x)=J_{\Lambda^{\prime\prime}}(y)$}\,.
$$
ii) Let $x=x_-ux_+\in A_2$ and set $y:=\overline{J_{\Lambda^{\prime\prime}}(x)}=y_-uy_+$. Then
\begin{equation}\label{equation}
\phi_\Lambda^{u,v}(x)=x_-vx_+\quad\text{and}\quad \phi_\Lambda^{u,v}(y_-uy_+)=y_-vy_+\,.
\end{equation}
Result ii) follows from $\lim_t y_-uy_+=x_-ux_+$ and $\lim_t y_-vy_+=x_-vx_+$.\\
iii) Since $f$ is uniformly continuous there exists $\delta>0$ such that
$$
d(x,y)\leq\delta\implies |f(x)-f(y)|\leq \varepsilon^\prime|\Lambda^\prime|^{-1}\,.
$$
Choose $t^\prime:=r+m^*$, with $\epsilon_{m^*}\leq\delta$ (see \eqref{distance}).
If  $x=x_-ux_+\in A_2(n)\backslash A_2(n|m)$, then \eqref{equation} also holds, and
$$
\psi_{f,r}(x)=\sum_{j\in \Lambda^\prime}\big(f(T^j y_-vy_+)-
f(T^j y_-uy_+)\big)\,.
$$
For all $j\in\Lambda^\prime$,
$$
d(T^jx_-ux_+,T^jy_-uy_+)\leq\delta \quad\text{and}\quad
d(T^jx_-vx_+,T^jy_-vy_+)\leq\delta\,,
$$
so that
$$
\big|\sum_{j\in \Lambda^\prime}\big(f(T^j\phi_{\Lambda}^{u,v}(x))-f(T^jx)\big)-\psi_{f,r}(x)\big|\leq2\varepsilon^\prime\,.
$$
iv) This follows from iii), taking into account \eqref{eqr} . 
\qed

\begin{lem}\label{lemequal}
For any interval $\Lambda$ and $u,v\in\Sigma_\Lambda$,
$$
\int_{A_2}{\rm e}^{\psi_f}\,d\nu=\nu(B_2)\,.
$$
\end{lem}

\prf Let $\Lambda=[k,\ell]$ and $[-m,m]\supset [k-n,\ell+n]$ with $n\geq t\geq t^\prime$ (see lemma \ref{lemJC}).
Set
$$
\underline{A}_2(n|m):=A_2(n)\backslash A_2(n|m)\quad\text{and}\quad \underline{B}_2(n|m):=B_2(n)\backslash B_2(n|m)\,.
$$
There is a bijection between $\underline{A}_2(n|m))$ and 
$\underline{B}_2(n|m)$ and
\begin{eqnarray*}
\mu_m(B_2(n))
&=& \mu_m(\underline{B}_2(n)|m))+\mu_m(B_2(n|m))\\
&=&\int_{\underline{A}_2(n|m)}{\rm e}^{\psi_f}\,d\mu_m+\mu_m(B_2(n|m))\,.
\end{eqnarray*}
From lemma \ref{lemJC},
$$
{\rm e}^{-(\varepsilon+2\varepsilon^\prime)}\int_{\underline{A}_2(n|m)}{\rm e}^{\psi_{f,r}}\,d\mu_m
\leq
\int_{\underline{A}_2(n|m)}{\rm e}^{\psi_{f}}\,d\mu_m
\leq
{\rm e}^{\varepsilon+2\varepsilon^\prime}\int_{\underline{A}_2(n|m)}{\rm e}^{\psi_{f,r}}\,d\mu_m\,.
$$
These results and lemma \ref{lemnotb}  are also true for $\Lambda+j=[k+j,\ell+j]$ as long as  $[-m,m]\supset  [k+j-n,\ell+j+n]$.
Therefore, taking into account  the continuity of $\psi_{f,r}$ and lemma \ref{lemnotb}, for fixed $r$,
$$
\lim_m\int_{\underline{A}_2(n|m)}{\rm e}^{\psi_{f,r}}\,d\nu_m=\lim_m\int_{A_2(n)}{\rm e}^{\psi_{f,r}}\,d\nu_m=\int_{A_2(n)}{\rm e}^{\psi_{f,r}}\,d\nu\,.
$$
The set $B_2(n)$ is a finite disjoint union of cylinder sets; therefore
$$
\lim_m\nu_m(B_2(n))=\nu(B_2(n))\leq{\rm e}^{\varepsilon+2\varepsilon^\prime}\int_{A_2(n)}{\rm e}^{\psi_{f,r}}\,d\nu
$$
and
$$
\lim_m\nu_m(B_2(n))=\nu(B_2(n))\geq{\rm e}^{-(\varepsilon+2\varepsilon^\prime)}\int_{A_2(n)}{\rm e}^{\psi_{f,r}}\,d\nu\,.
$$
By monotonicity, taking  the limit $n\ra\infty$, 
$$
{\rm e}^{-(\varepsilon+2\varepsilon^\prime)}\int_{A_2}{\rm e}^{\psi_{f,r}}\,d\nu\leq
\nu(B_2)\leq {\rm e}^{\varepsilon+2\varepsilon^\prime}\int_{A_2}{\rm e}^{\psi_{f,r}}\,d\nu\,.
$$
Taking the limit $t\ra\infty$ (see lemma \ref{lemJC} ii), and then the limit $r\ra\infty$, one gets
$$
\nu(B_2)=\int_{A_2}{\rm e}^{\psi_f}\,d\nu
$$
since $\varepsilon$ and $\varepsilon^\prime$ are arbitrary.
\qed

\begin{thm}\label{thmconformal}
For any $\beta >1$, if $f\in \B(\Sigma^\beta)$ verifies Bowen's condition, then the unique equilibrium measure $\nu$ is a Gibbs conformal measure.
\end{thm}

\prf
Since $\psi_f$ is trivial on $A_1$,
$$
\nu\big(\phi_\Lambda^{u,v}([u])\big)=\nu(A_1)+\nu(B_2)=\nu(A_1)+\int_{A_2}{\rm e}^{\psi_f}\,d\nu=\int_{[u]}{\rm e}^{\psi_f}\,d\nu\,.
$$
\qed

\end{document}